\documentclass[english]{elsarticle}
\usepackage{mathptmx}
\usepackage[T1]{fontenc}
\usepackage[latin1]{inputenc}
\usepackage{fancyhdr}
\usepackage{amssymb}
\usepackage{graphicx}

\makeatletter

\providecommand{\tabularnewline}{\\}

\@ifundefined{date}{}{\date{}}
\usepackage{rotating}
\usepackage{amssymb}
\usepackage{amsthm}
\usepackage{booktabs}

\usepackage{fancyheadings} 

\@ifundefined{showcaptionsetup}{}{%
 \PassOptionsToPackage{caption=false}{subfig}}
\usepackage{subfig}
\makeatother

\usepackage{babel}
\begin{document}
\fancyhead{} \chead{PREPRINT: DO NOT CITE OR DISTRIBUTE}  \renewcommand{\headrulewidth}{0pt}

\title{Multiwavelet density estimation\thispagestyle{fancy}}

\author[rvt]{Judson B. Locke\corref{cor1}}

\ead{lockej@my.fit.edu}

\author[focal]{Adrian M. Peter\corref{cor1}}

\ead{apeter@fit.edu}

\cortext[cor1]{Corresponding authors. Send correspondence to .}

\address[rvt]{Florida Institute of Technology, Department of Physics and Space
Sciences, 150 West University Blvd., Melbourne, FL 32901, USA}

\address[focal]{Florida Institute of Technology, Department of Engineering Systems,
150 West University Blvd., Melbourne, FL 32901, USA}
\begin{abstract}
Accurate density estimation methodologies play an integral role in
a variety of scientific disciplines, with applications including simulation
models, decision support tools, and exploratory data analysis. In
the past, histograms and kernel density estimators have been the predominant
tools of choice, primarily due to their ease of use and mathematical
simplicity. More recently, the use of wavelets for density estimation
has gained in popularity due to their ability to approximate a large
class of functions, including those with localized, abrupt variations.
However, a well-known attribute of wavelet bases is that they can
not be simultaneously symmetric, orthogonal, and compactly supported.
Multiwavelets---a more general, vector-valued, construction of wavelets---overcome
this disadvantage, making them natural choices for estimating density
functions, many of which exhibit local symmetries around features
such as a mode. We extend the methodology of wavelet density estimation
to use multiwavelet bases and illustrate several empirical results
where multiwavelet estimators outperform their wavelet counterparts
at coarser resolution levels.\end{abstract}
\begin{keyword}
Orthogonal series density estimation, non-parametric density estimation,
wavelets, multiwavelets
\end{keyword}
\maketitle

\section{Introduction\label{sec:Introduction}}

Many applications require estimating the underlying probability density
function (PDF) of a finite sample making minimal or no assumptions
about the generating function. Having an accurate model of the underlying
PDF enables one to understand the structure of the data and carry
out more advanced statistical analysis such as classification, confidence
modeling, and clustering. Here we introduce for the first time a new
class of density estimators based on a series expansion of multiwavelets
\citep{Alpert1993,Goodman1994}. Multiwavelets can be created to retain
all of the desirable properties of regular wavelets and, in addition,
exhibit very desirable properties which wavelets do not: simultaneous
symmetry, compact support, and orthogonality. These properties make
multiwavelet density estimation (MWDE) well-suited for reconstructing
a wide class of density functions, particularly those that exhibit
local or global symmetries. 

The primary contribution of our work is to introduce for the \emph{first
time} the use of multiwavelets for density estimation. We also empirically
compare MWDE performance to regular wavelet density estimation (WDE).
Our focus is not to pit multiwavelets versus wavelets for the task
of density estimation. We are still in the early stages of our research
and continue to explore pros and cons associated with multiwavelets
used as bases for density estimation. For example, the methodology
we present here is strictly a \emph{linear} multiwavelet estimator;
hence, we do not discuss issues of thresholding the multiwavelet bases.
Implementing an effective thresholding technique will yield sparser
representations and should make MWDE more directly comparable to WDE.
The richer mathematical properties afforded by multiwavelets demand
we investigate their use for important applications such as density
estimation.

\subsection{Relevant Work}

There are many well-studied density estimation techniques which we
can loosely categorized into the taxonomy of parametric, semi-parametric,
and nonparametric estimators. Of these, nonparametric models are the
most data-driven, requiring little to no assumptions about the underlying
generative model of the data. These models include the ubiquitous
histogram described by \citet{Silverman1998} and the well-established
kernel density estimators \citep{Scott01}. Though introduced in the
1960s by \citet{Cencov62} and later described by \citet{Watson1969}
and improved by \citet{Anderson1978}, orthogonal series estimation
(OSE) lagged in popularity due to lack of suitable bases for the series
expansion. Despite this lag, work was done on OSE by \citet{Hall1986}
and \citet{Ahmad1982} to investigate the convergence rate and integrated
mean square error properties, respectively. Until about 25 years ago,
the series expansions used for OSE were essentially limited to Fourier
bases (i.e. sines and cosines) \citep{Kronmal1968} or orthogonal
polynomials, e.g. Hermite \citep{Schwartz1967} and Laguerre \citep{Izenman1991}.
The main downfall of these bases is their infinite support, demanding
a large number of terms in the series expansion to accurately approximate
complex densities containing multiple modes and abrupt variations.
With the advent and growing use of wavelets, we are now seeing more
uses of OSE \citep{Garcia2012,Provost2012,Cheng2009,Caudle09}. In
fact, \citet{Peter08} show wavelet density estimators (WDE) often
outperform many other nonparametric density estimators. The main advantage
comes from the fact wavelets can be constructed with the convenient
property of compact support \citep{Daubechies92}, allowing them to
easily and accurately represent functions with discontinuities and
other abrupt local phenomena. In addition, WDE can be extended to
non-linear estimators through the use of wavelet function coefficient
thresholding introduced by \citet{Donoho96}, which allows us to represent
the density with a sparse set of coefficients while retaining accuracy.

Unfortunately, wavelets can not be simultaneously compact, symmetric,
and orthogonal \citep{Daubechies92}. Many common analytic densities
exhibit some form of symmetry either globally (e.g. Gaussian and Laplacian)
or locally about their modes. The immediate consequence of wavelets
lacking symmetry (if they also want to be orthogonal and compact)
is that more of them will be required, either via a multiresolution
expansion or a very fine single level expansion, to reproduce the
symmetries in the underlying densities. Doing away with compactness
and orthogonality allows us to have wavelet bases that are symmetric,
but these relinquished features are the very properties that made
wavelets an attractive choice over trigonometric or polynomial bases.
This leaves us with only one choice of bases satisfying all these
desirable properties: multiwavelets. 

Multiwavelets are a more general, vector-valued construction of wavelets
\citep{Alpert1993,Goodman1993,Goodman1994}. When used in a series
expansion, they utilize multiple basis functions at every translate
and at each resolution level. Like wavelets, multiwavelets can be
constructed to have compact support, allowing them to faithfully model
local discontinuities. Furthermore, they can be orthogonal, making
coefficient estimation simple. Unlike wavelets, they can be symmetric
while maintaining their compactness and orthogonality, allowing them
to model local and global symmetries at coarser resolutions than wavelets.
In this paper, we extend the concept of WDE to multiwavelet density
estimation and demonstrate the utility of multiwavelets to model a
variety of distributions. To our knowledge, this is the first use
of multiwavelets for density estimation.

The remainder of this paper is organized as follows. The next section
provides background knowledge of wavelets, WDE, and multiwavelets.
In \S\ref{sec:MWDE}, we show how to construct a density estimator
with multiwavelets. In \S\ref{sec:Experiments}, we demonstrate the
capabilities of MWDE using a wide variety of multiwavelet families
and compare MWDE to the conventional WDE. Finally, we conclude with
some observations and directions for future research.

\section{Review of Wavelets and Multiwavelets\label{sec:Theory}}

\subsection{Standard Wavelets and Multiresolution Analysis\label{subsec:Wavelets}}

Wavelets are refinable functions whose values are solutions to the
dilation equation
\begin{equation}
\phi\left(x\right)=\sqrt{m}\sum_{k\in\mathbb{Z}}h_{k}\phi\left(mx-k\right),\label{eq:WaveletDilationEquation}
\end{equation}
where $\phi\left(x\right)$ is called the father wavelet (a.k.a. scaling
function), $m\in\mathbb{Z}$ is the dilation factor (in our case,
and in most practical cases, $m=2$), $h_{k}\in\mathbb{R}$ are low-pass
coefficients (the ``filter'') defining the wavelet, and $k$ are
integer translates of the father wavelet across the domain. The dilation
factor is further controlled by a choice of resolution level, $j\in\mathbb{Z}$,
that expands or contracts the wavelet. Hence we will typically denote
the normalized basis function at resolution $j$ and translate $k$
as $\phi_{j,k}=2^{j/2}\phi(2^{j}x-k)$. At a chosen resolution level
$j$, the father wavelet and its integer translates form a basis for
the function space $V_{j}$ which is a subspace of the space of all
square integrable functions $L^{2}(\mathbb{R})$. These father wavelets
capture the ``smooth'' or ``averaging'' properties of functions.
Correspondingly, one can construct a set of mother wavelets which
model the details (i.e. oscillating properties) of functions. These
form a basis for the space $W{}_{j}$, which consists of functions
orthogonal to $V_{j}$. The mother wavelets $\psi\left(x\right)$
(a.k.a. wavelet functions) are the members of $W$ and are found using
the high-pass filter coefficients $g_{k}$ and by solving the wavelet,
two-scale equation
\begin{equation}
\psi\left(x\right)=\sqrt{2}\sum_{k\in\mathbb{Z}}g_{k}\phi\left(2x-k\right).\label{eq:WaveletEquation}
\end{equation}
Eqs. \ref{eq:WaveletDilationEquation} and \ref{eq:WaveletEquation}
only provide values for $\phi$ and $\psi$ evaluated on the dyadic
rationals. If we want to evaluate $\phi\left(x\right)$ at any $x\in\mathbb{R}$,
then we simply interpolate using known values of $\phi$ on the dyadic
rationals. When this arises in the context of density estimation,
we use a cubic spline interpolant.

There are several families of wavelets---Daubechies, Coiflets, Symlets,
just to name a few---all of which are constructed from and defined
by unique sets of coefficients. Wavelets can be constructed to have
the convenient properties of orthogonality across their integer translates
and compact support; however, these properties can not be had simultaneously
with perfect symmetry \citep{Daubechies92}. As discussed in \S\ref{subsec:Multiwavelets},
this drawback can be addressed through the use of multiwavelets. 

Scaling and wavelet bases can be brought together to represent functions
in a multiresolution expansion. Given a function $f\in L^{2}\left(\mathbb{R}\right)$,
a multiresolution analysis (MRA) of that function yields 

\begin{equation}
f\left(x\right)=\sum_{k\in\mathbb{Z}}\alpha_{j_{0},k}\phi_{j_{0},k}\left(x\right)+\sum_{j=j_{0}}^{\infty}\sum_{k\in\mathbb{Z}}\beta_{j,k}\psi_{j,k}\left(x\right),\label{eq:WaveletDensityExpansion-1}
\end{equation}
where $\alpha_{j_{0},k}\in\mathbb{R}$ are the father wavelet coefficients,
and $\beta_{j,k}\in\mathbb{R}$ are the mother wavelet coefficients.
The lowest resolution level used for function approximation is $j_{0}\in\mathbb{Z}$,
and all other resolutions are $j\in\mathbb{Z}$ subject to $j\geq j_{0}$.
In reality, though, only a finite set of resolutions will be utilized
such that $j_{0}\leq j\leq J<\infty$. There are several approaches
that can be applied for choosing $j_{0}$ and $J$; in the context
of density estimation, we refer interested readers to \citet{Vidakovic99}.
The MRA in $L^{2}(\mathbb{R})$ is a doubly-infinite sequence of nested
subspaces $V_{j\in\mathbb{Z}}$
\begin{equation}
\cdots\subset V_{-2}\subset V_{-1}\subset V_{0}\subset V_{1}\subset V_{2}\subset\cdots
\end{equation}
such that $\bigcap_{j}V_{j}=\left\{ 0\right\} $ and $\overline{\bigcup_{j}V_{j}}=L^{2}(\mathbb{R})$,
allowing $V_{j}$ to be used as a basis. The scaling and wavelet bases
relationship is such that $W_{j}$ is orthogonal to $V_{j}$, and
at any resolution level $j$ we have
\begin{equation}
V_{j}=V_{j-1}\oplus W_{j-1},
\end{equation}
justifying the expansion in Eq. \ref{eq:WaveletDensityExpansion-1}.

\subsection{Wavelet Density Estimation\label{subsec:WDE}}

Wavelet density estimation is a specific incarnation of orthogonal
series density estimation (OSE), where one expands the density function
$p(x)$ in a multiresolution wavelet basis:

\begin{equation}
p\left(x\right)=\sum_{k\in\mathbb{Z}}\alpha_{j_{0},k}\phi_{j_{0},k}\left(x\right)+\sum_{j=j_{0}}^{\infty}\sum_{k\in\mathbb{Z}}\beta_{j,k}\psi_{j,k}\left(x\right).\label{eq:WaveletDensityExpansion}
\end{equation}
The objective now becomes determining the coefficients $\alpha_{j_{0},k}$
and $\beta_{j,k}$ from a given independent and identically distributed
(i.i.d.) sample $X=\left\{ X_{i}\right\} _{i=1}^{N}$. The standard
approach---though there are others \citep{Peter08}---is to simply
project the density function onto the basis expansion:
\begin{equation}
\alpha_{j_{0},k}=\left\langle p,\phi_{j_{0},k}\right\rangle =\int p\left(x\right)\phi_{j_{0},k}\left(x\right)dx.\label{eq:WaveletInnerProduct}
\end{equation}
This allows us to interpret the inner product of Eq. \ref{eq:WaveletInnerProduct}
as an expectation to find
\begin{equation}
\alpha_{j_{0},k}=\int p\left(x\right)\phi_{j_{0},k}\left(x\right)dx=\mathcal{E}\left[\phi_{j_{0},k}\left(x\right)\right],
\end{equation}
where $\mathcal{E}\left[\cdot\right]$ is the expectation operator.
Assuming a finite sample, the expectation is approximated by the sample
mean. Hence, the coefficients $\alpha_{j_{0},k}$ are estimated by
\begin{equation}
\hat{\alpha}_{j_{0},k}=\frac{1}{N}\sum_{i=1}^{N}\phi_{j_{0},k}\left(X_{i}\right),
\end{equation}
and similarly for the wavelet function coefficients $\beta_{j,k}$
as 
\begin{equation}
\hat{\beta}_{j,k}=\frac{1}{N}\sum_{i=1}^{N}\psi_{j,k}\left(X_{i}\right),
\end{equation}
to arrive at the approximation $\hat{p}\left(x\right)$ given by
\begin{equation}
\hat{p}\left(x\right)=\sum_{k\in\mathbb{Z}}\hat{\alpha}_{j_{0},k}\phi_{j_{0},k}\left(x\right)+\sum_{j=j_{0}}^{J}\sum_{k\in\mathbb{Z}}\hat{\beta}_{j,k}\psi_{j,k}\left(x\right).
\end{equation}

There are no guarantees, however, that the resulting density estimate
$\hat{p}$ will satisfy the necessary properties of a density function
(namely, $\int\hat{p}=1$ and $\hat{p}\geq0$). So, once the density
has been estimated, a post-processing normalization procedure \citep{Gajek1986}
is typically performed to achieve these properties. It is worth noting
there are ways to guarantee the resulting $\hat{p}$ is positive,
such as using positive wavelets as in \citet{Walter1999}, and even
ways to ensure the resulting density will be positive and integrate
to one, namely by estimating not $p$, but $\left(\sqrt{p}\right)^{2}=p$
with certain restrictions on the coefficients \citep{Pinheiro97,Vidakovic99,Peter08}.

The utility of representing a density function using an MRA stems
from the ability to threshold detail coefficients and gain an economical
(in the sense of sparse coefficients), yet accurate estimator. These
coefficients can be set to zero via ``hard thresholding'' techniques.
Similarly, the larger mother wavelet function coefficients can be
shrunk toward zero to reduce their contribution to the reconstruction,
making the resulting estimate smoother; this is generally called ``soft
thresholding.'' In \citet{Donoho96} and \citet{Donoho98}, thresholding
and its implications on convergence are analyzed to provably show
WDE to be optimal under mini-max criteria. When coefficient thresholding
is employed, WDE is a class of nonlinear density estimators. Our present
focus is to introduce the use of multiwavelets in density estimation.
To this end, we do not address the issue of thresholding coefficients
when using multiwavelets as the bases of our density estimator.

\subsection{Multiwavelets\label{subsec:Multiwavelets}}

Multiwavelets are a more general, vector-valued constructions of wavelets.
Excellent foundations of multiwavelet theory are given in \citet{Strela1996}
and \citet{Keinert04}, with origins traced back to \citet{Alpert1993},
\citet{Goodman1993}, and \citet{Goodman1994}. A multiscaling function
$\underline{\phi}\left(x\right)$ is a vector-valued refinable function
with multiplicity $r\in\mathbb{Z}^{+}$ of the form
\begin{equation}
\underline{\phi}\left(x\right)=\left(\begin{array}{c}
\phi_{1}\left(x\right)\\
\vdots\\
\phi_{r}\left(x\right)
\end{array}\right),
\end{equation}
and satisfying the refinement equation
\begin{equation}
\underline{\phi}\left(x\right)=\sqrt{m}\sum_{k\in\mathbb{Z}}H_{k}\underline{\phi}\left(mx-k\right),\hspace{1em}k,m\in\mathbb{Z},\hspace{1em}m\geq2,
\end{equation}
where $H_{k}$ are low-pass $r\times r$ matrices called the ``recursion
coefficients'' defining the multiscaling function and collectively
referred to as the ``multifilter,'' paralleling the convention in
standard wavelets. As before, we are concerned only with the dyadic
case: $m=r=2$. In fact, all density estimations with multiwavelets
presented in this paper are performed, for the sake of simplicity,
with multiwavelets of multiplicity $r=2$. Multiscaling functions
can be generated using a cascade algorithm of the same form as the
standard wavelet cascade algorithm, but with matrix coefficients:
\begin{equation}
\underline{\phi}^{\left(n\right)}\left(x\right)=\sqrt{2}\sum_{k\in\mathbb{Z}}H_{k}\underline{\phi}^{\left(n-1\right)}\left(2x-k\right),
\end{equation}
with an orthogonal $\underline{\phi}^{\left(0\right)}$, such as the
Haar multiscaling function, paralleling the conventional Haar scaling
function.

Mother multiwavelet functions $\underline{\psi}\left(x\right)$ (in
this case, of multiplicity 2) can be created using the multiwavelet
equation (paralleling Eq. \ref{eq:WaveletEquation} for standard wavelets):
\begin{equation}
\underline{\psi}\left(x\right)=\sqrt{2}\sum_{k\in\mathbb{Z}}G_{k}\underline{\phi}\left(2x-k\right),\label{eq:MultiwaveletEquation}
\end{equation}
where $G_{k}$ are high-pass $r\times r$ matrices.

Analogous to wavelets in \S\ref{subsec:WDE}, MRAs can be constructed
for multiwavelets using the same procedure yielding the set $\left\{ \underline{\psi}_{j,k}:j,k\in\mathbb{Z}\right\} $.
This set is a basis for the space $W$. Like with traditional scalar
wavelets, the multiscaling functions form a basis for the spaces $V_{j\in\mathbb{Z}}$,
which are orthogonal to $W_{j\in\mathbb{Z}}$. As previously stated,
we are not presently concerned with coefficient thresholding, so all
multiwavelet density estimations in this paper are made using strictly
multiscaling functions.

\begin{figure}
\noindent \begin{centering}
\subfloat[]{\begin{centering}
\includegraphics[width=2in]{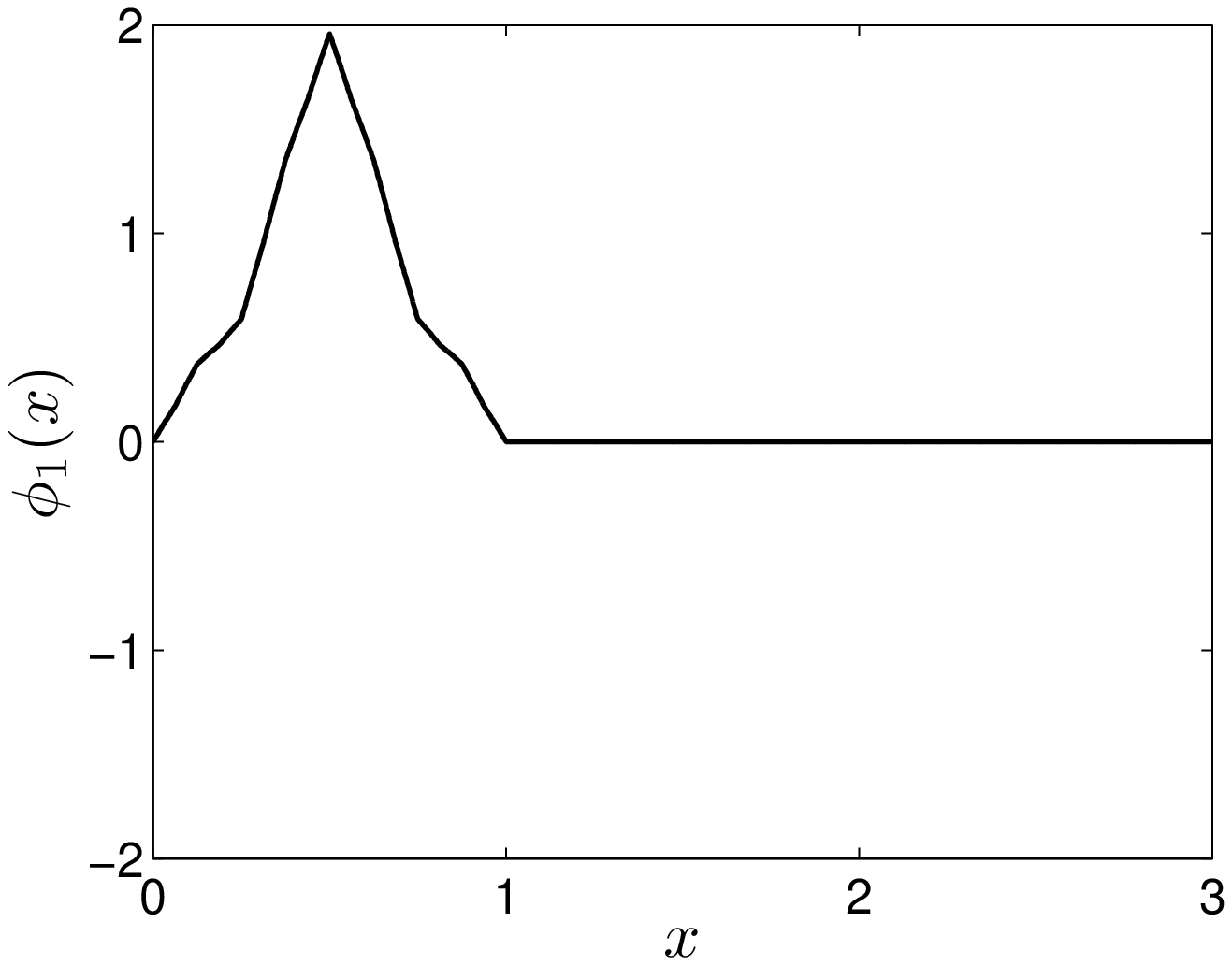}
\par\end{centering}

}\subfloat[]{\begin{centering}
\includegraphics[width=2in]{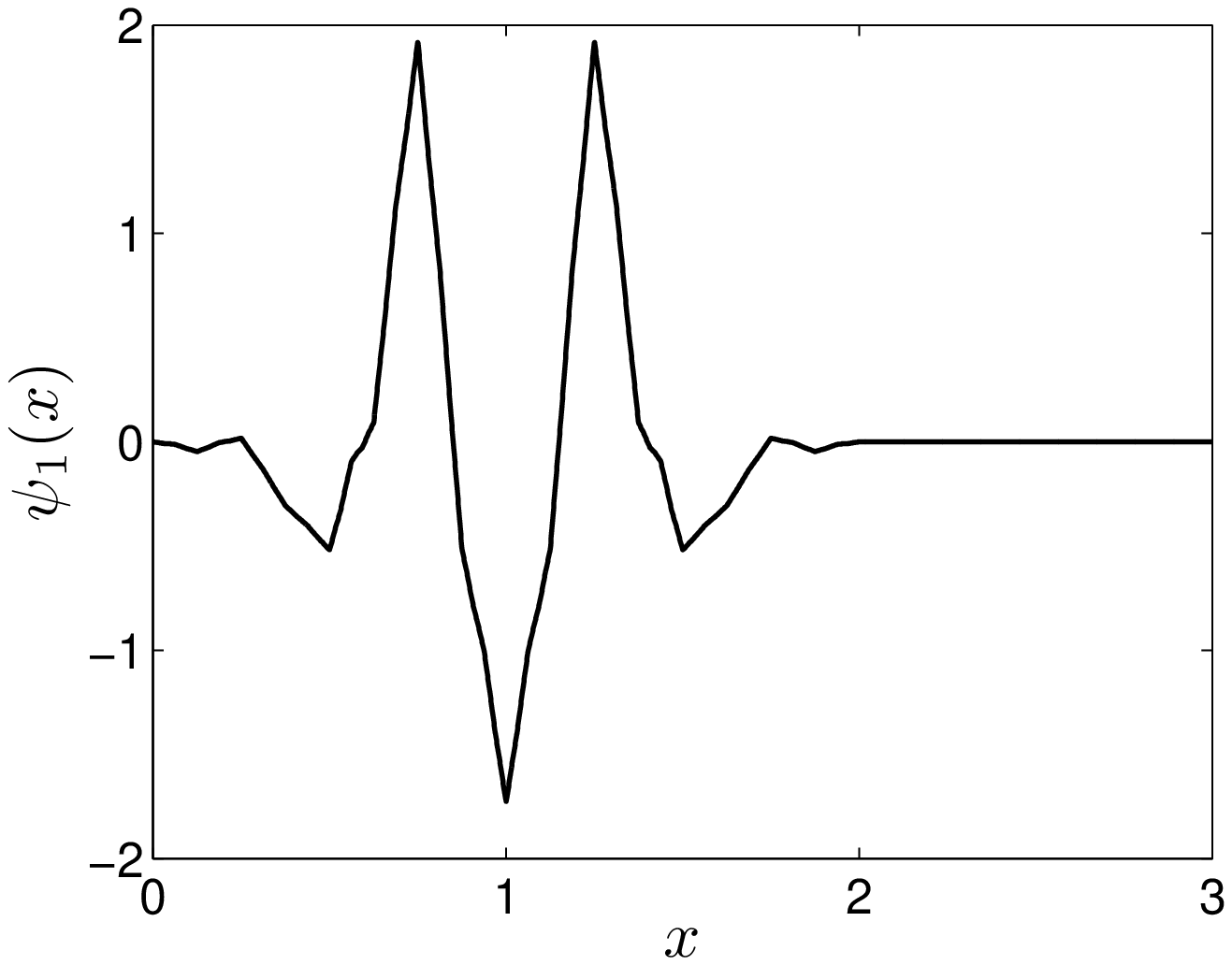}
\par\end{centering}

}
\par\end{centering}

\noindent \begin{centering}
\subfloat[]{\begin{centering}
\includegraphics[width=2in]{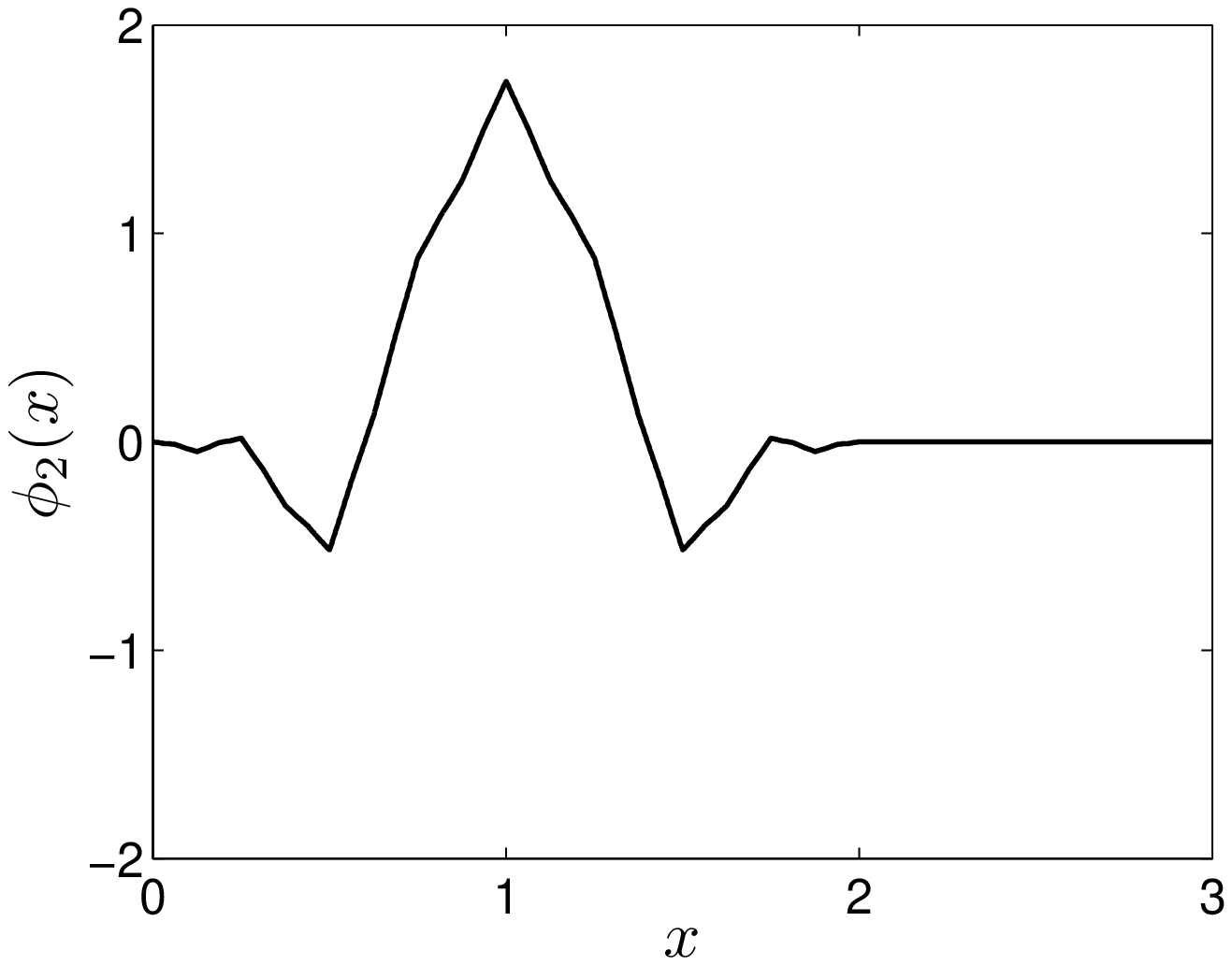}
\par\end{centering}

}\subfloat[]{\begin{centering}
\includegraphics[width=2in]{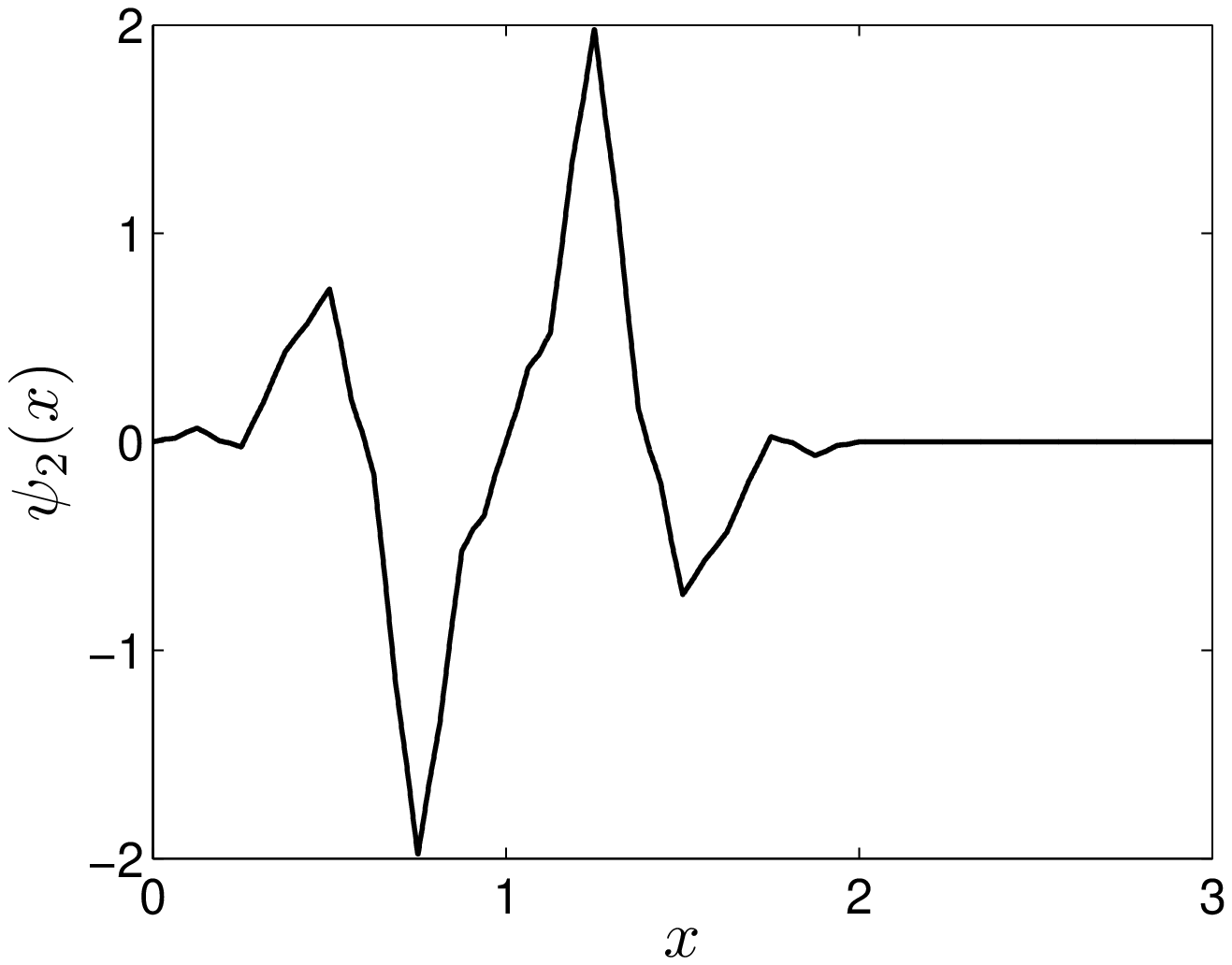}
\par\end{centering}

}
\par\end{centering}

\noindent \centering{}\caption{The DGHM multiwavelet has symmetric multiscaling functions \citet{DGHM1996}
and symmetric/antisymmetric multiwavelet functions \citet{StrangStrela1995}.
\label{fig:DGHM}}
\end{figure}

\begin{figure}
\noindent \begin{centering}
\subfloat[]{\noindent \begin{centering}
\includegraphics[width=2in]{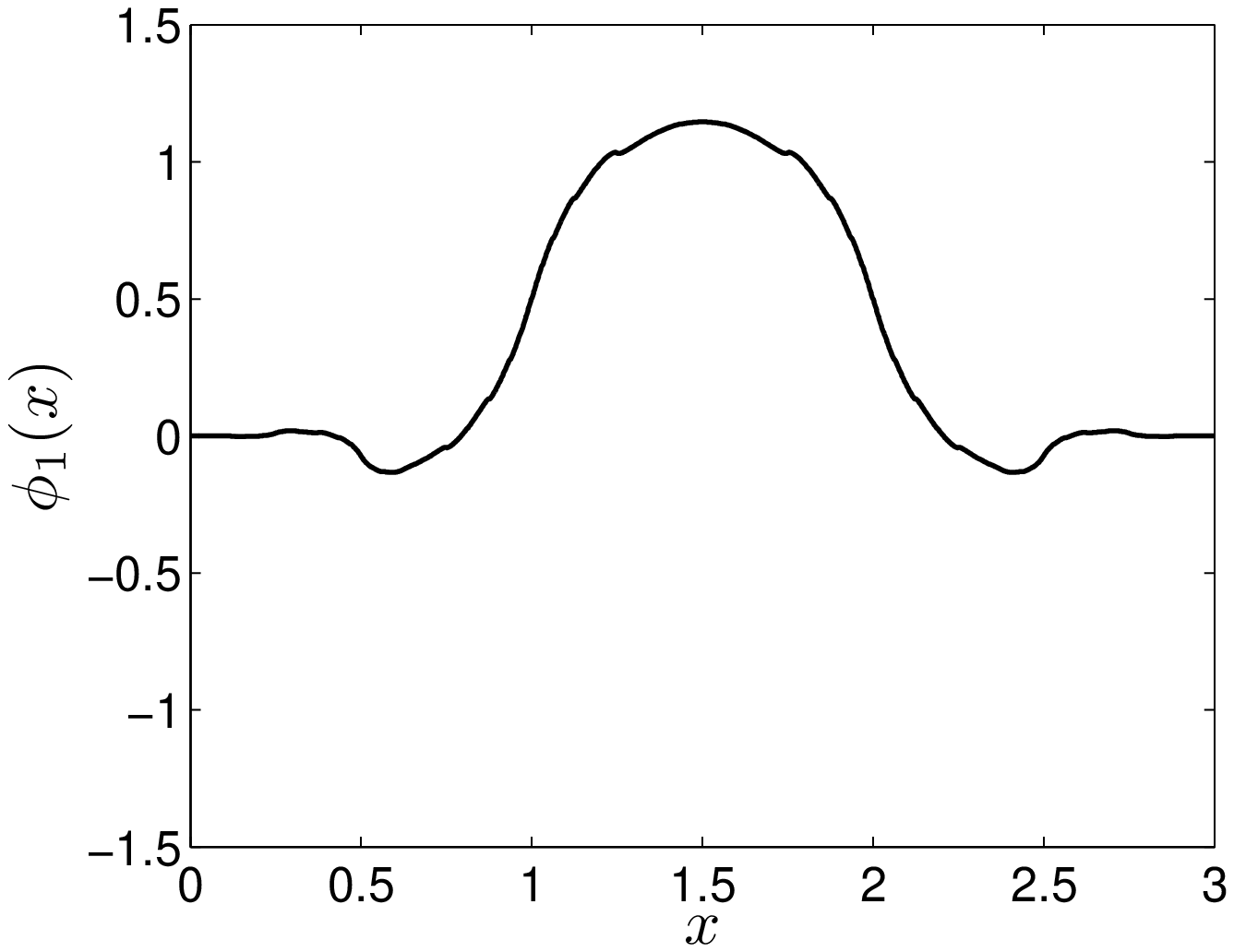}
\par\end{centering}

}\subfloat[]{\noindent \begin{centering}
\includegraphics[width=2in]{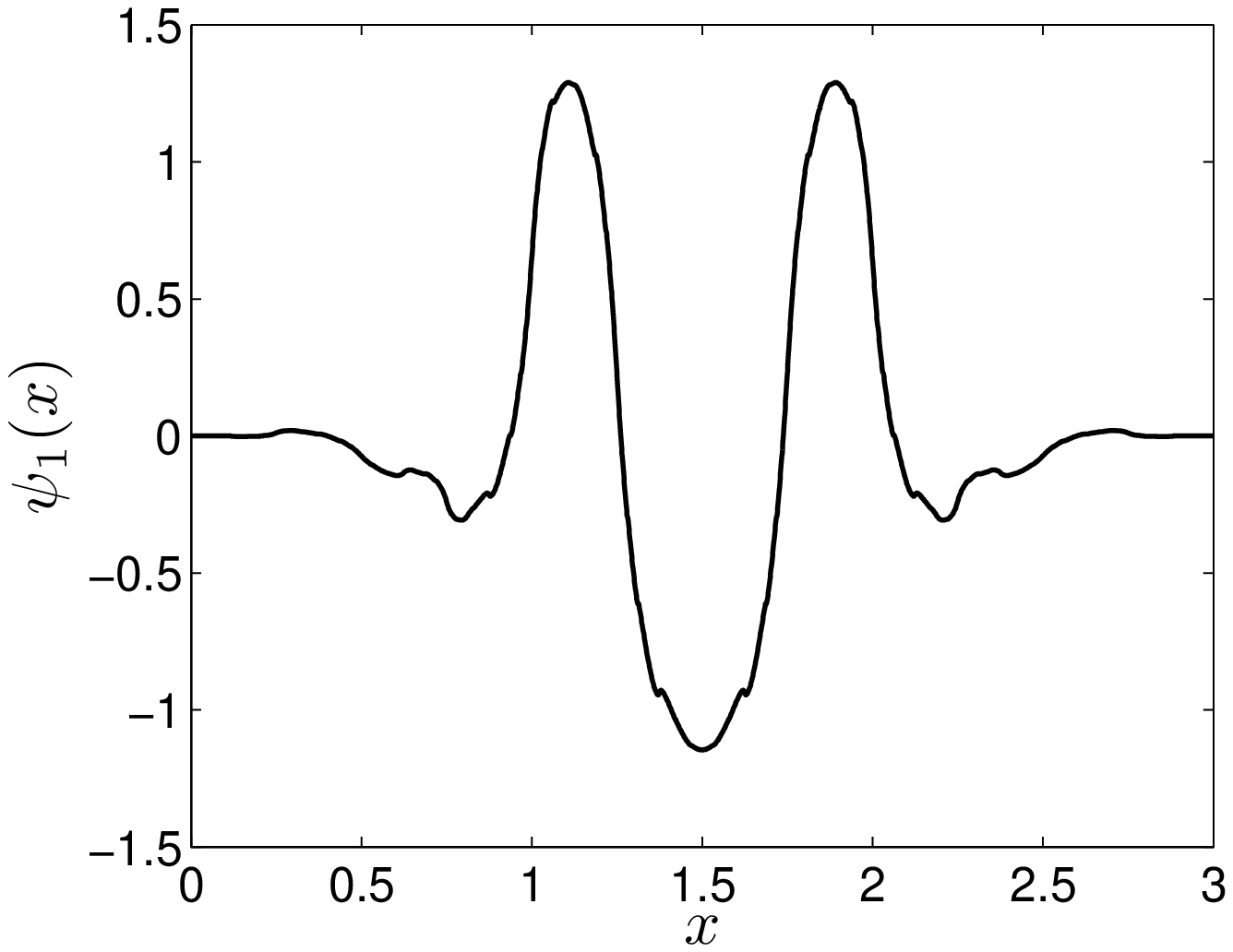}
\par\end{centering}

}
\par\end{centering}

\noindent \begin{centering}
\subfloat[]{\noindent \begin{centering}
\includegraphics[width=2in]{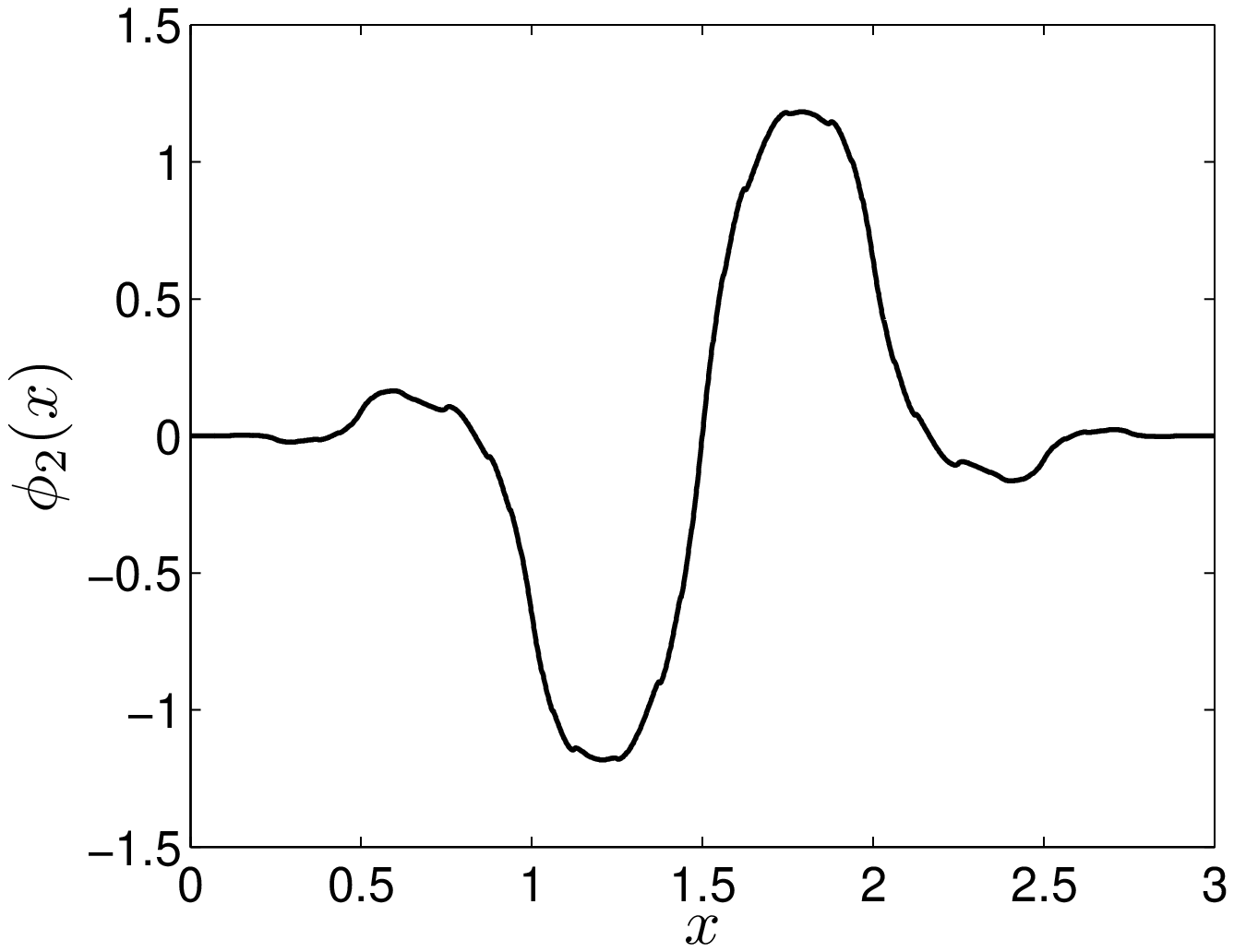}
\par\end{centering}

}\subfloat[]{\noindent \begin{centering}
\includegraphics[width=2in]{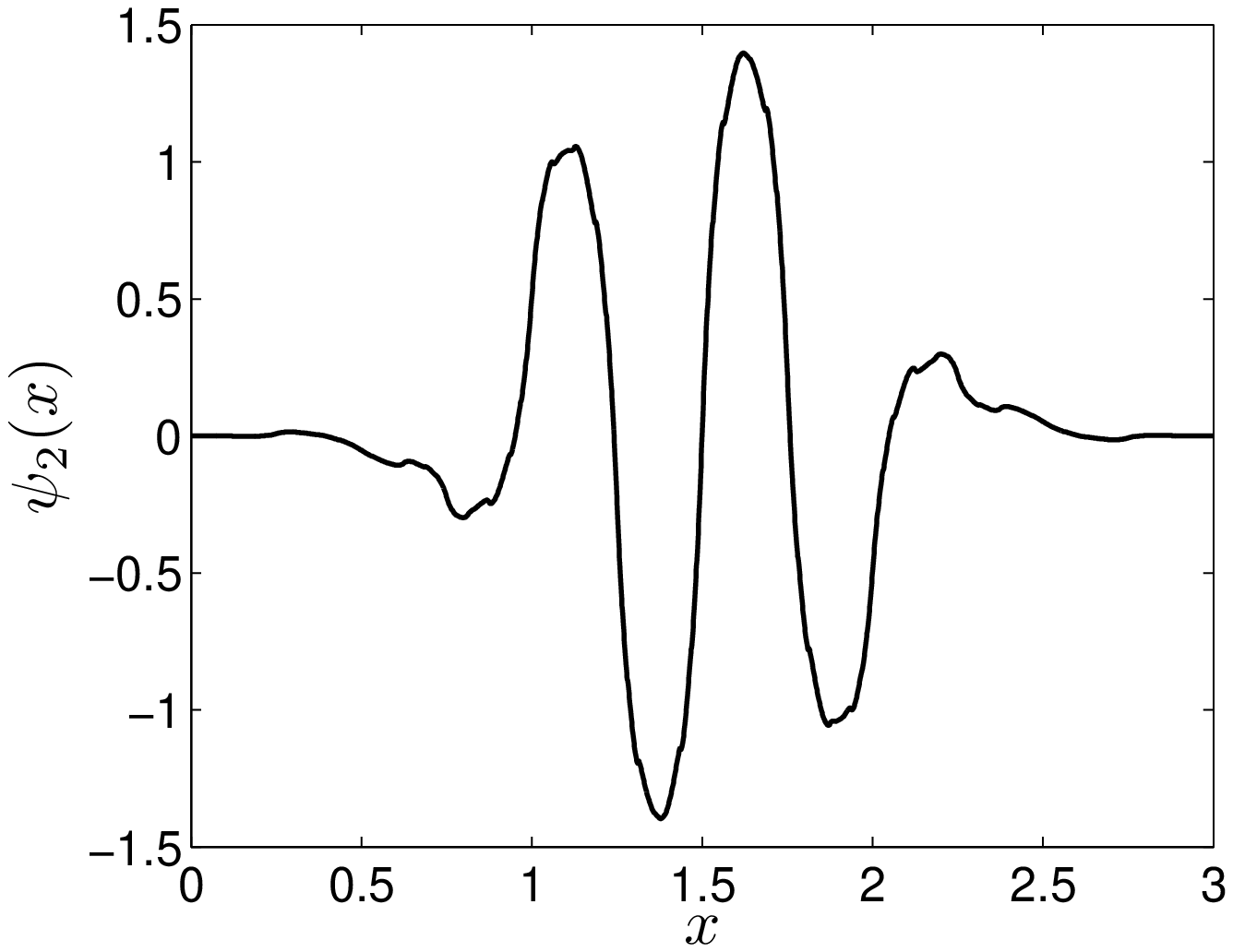}
\par\end{centering}

}
\par\end{centering}

\noindent \centering{}\caption{The STT multiwavelet \citet{STT2000} has symmetric/antisymmetric
multiscaling and multiwavelet functions. \label{fig:STT}}
\end{figure}

\begin{figure}
\begin{centering}
\begin{tabular}{cc}
\subfloat[]{\begin{centering}
\includegraphics[width=2in]{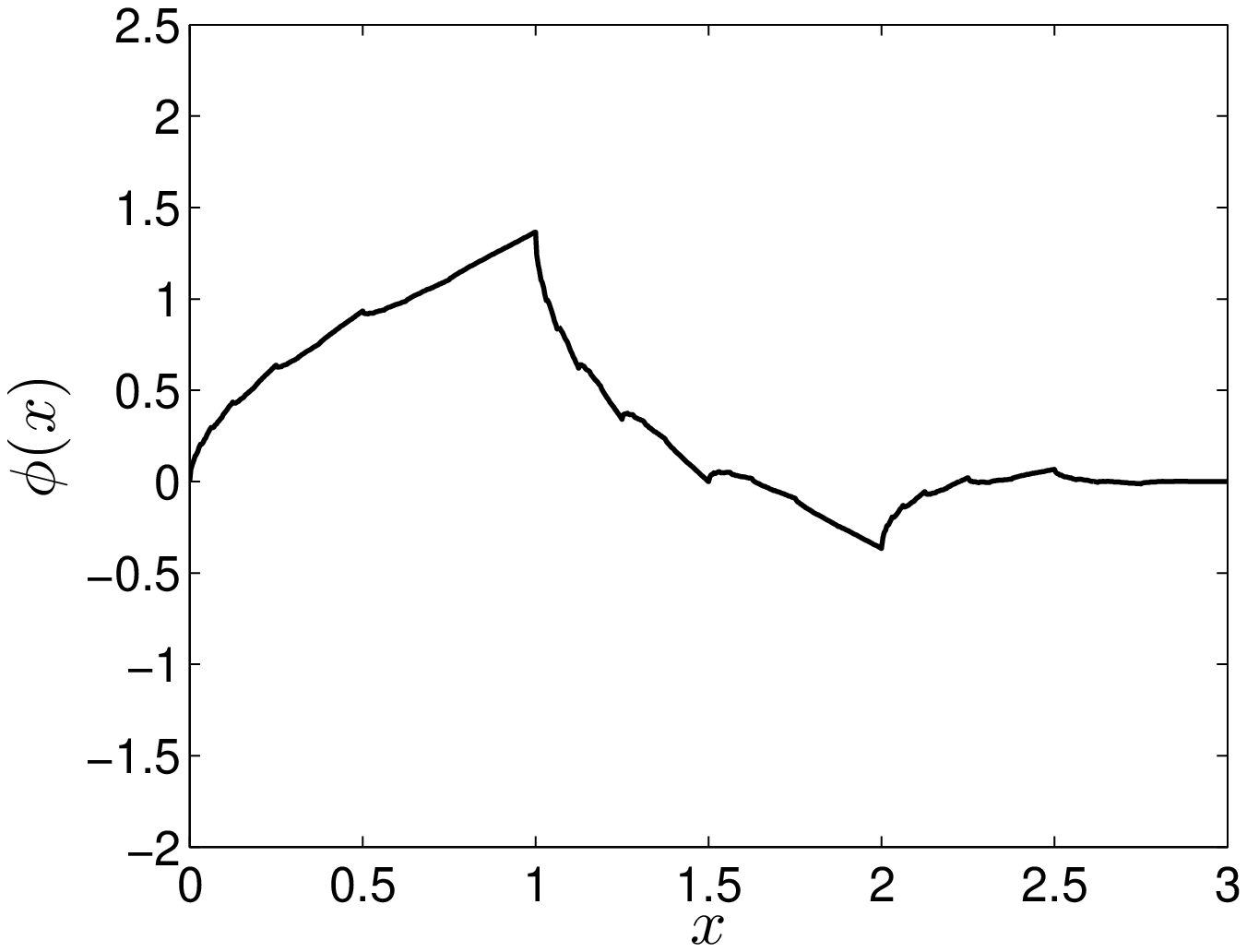}
\par\end{centering}

} & \subfloat[]{\begin{centering}
\includegraphics[width=2in]{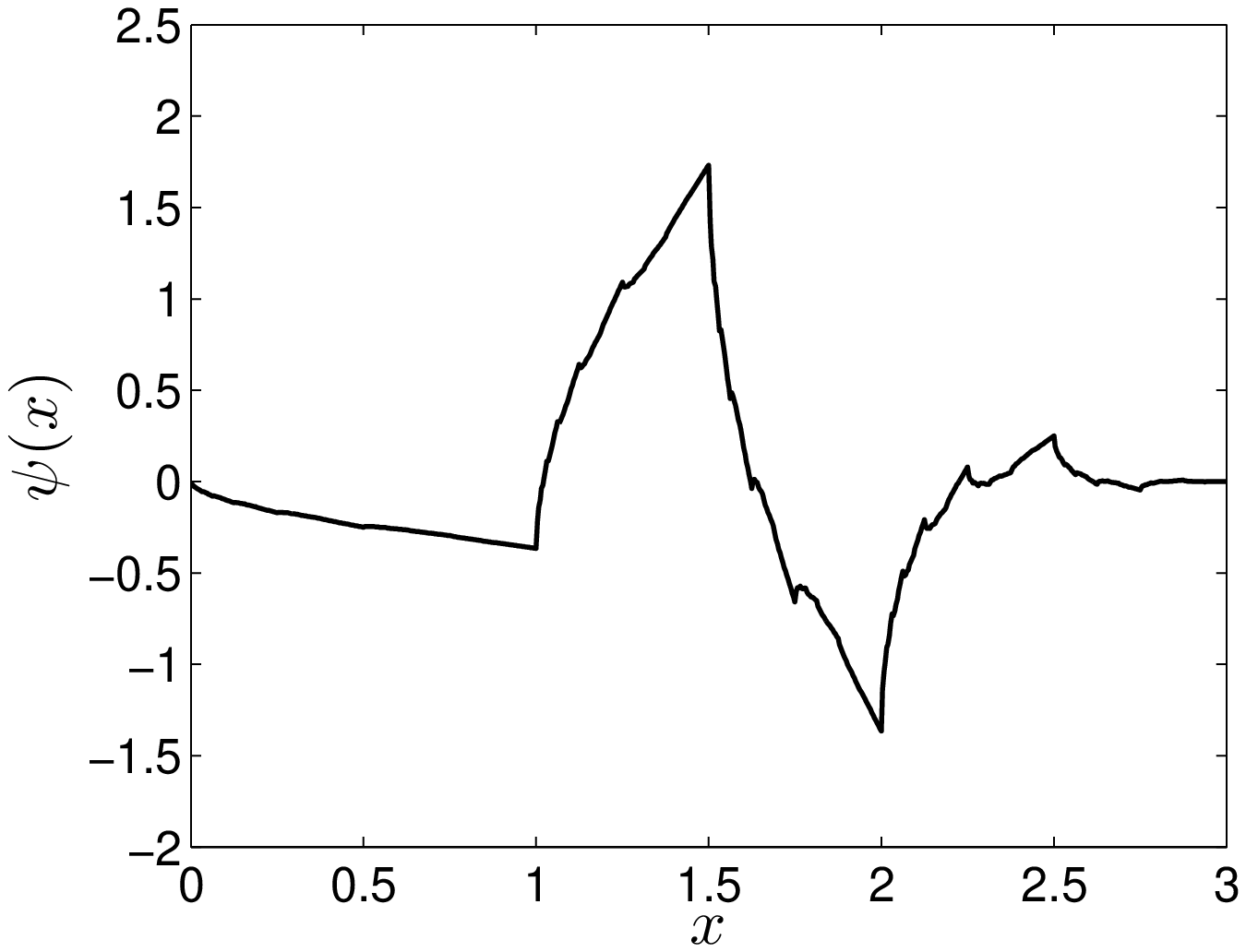}
\par\end{centering}

}\tabularnewline
\subfloat[]{\centering{}\includegraphics[width=2in]{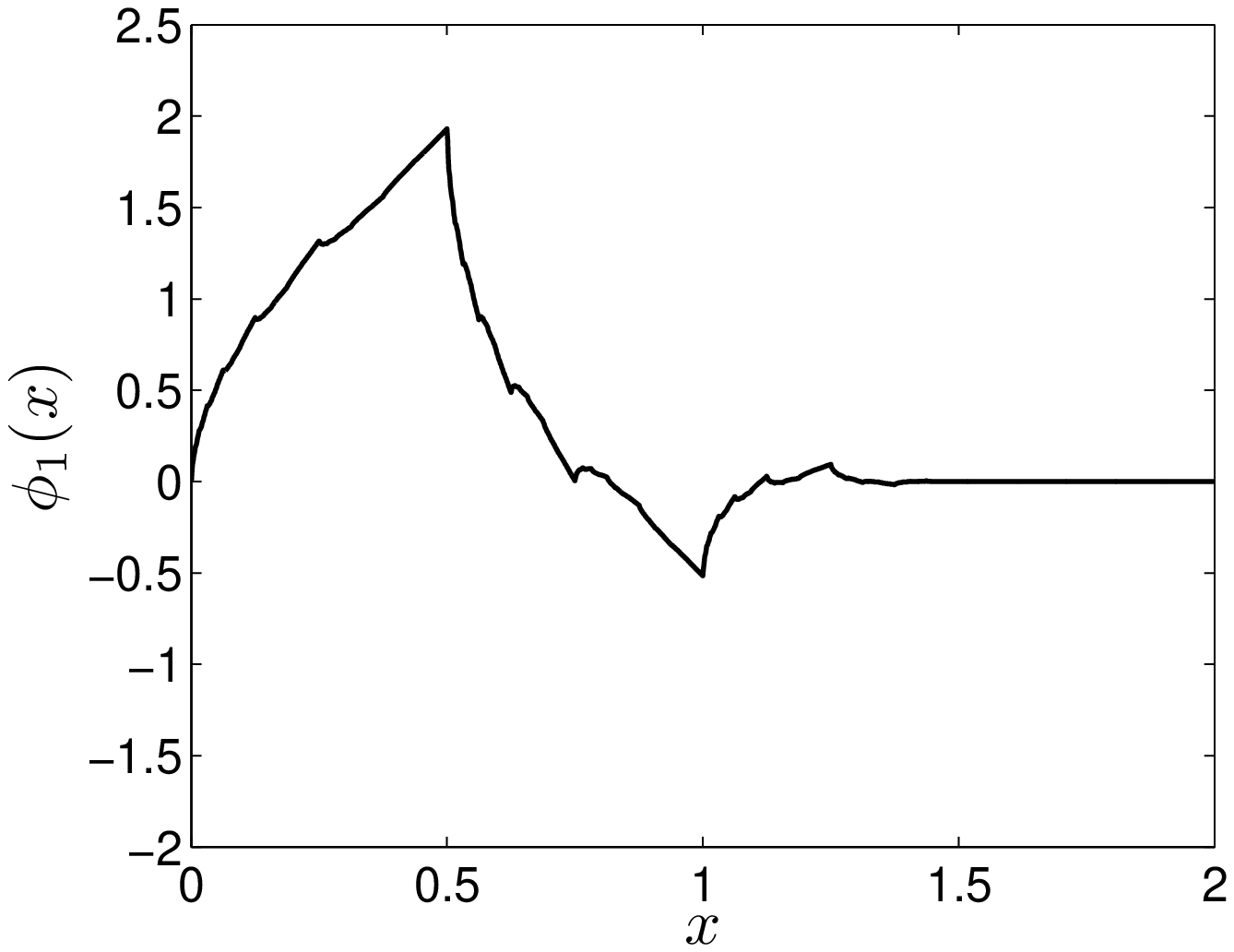}} & \subfloat[]{\begin{centering}
\includegraphics[width=2in]{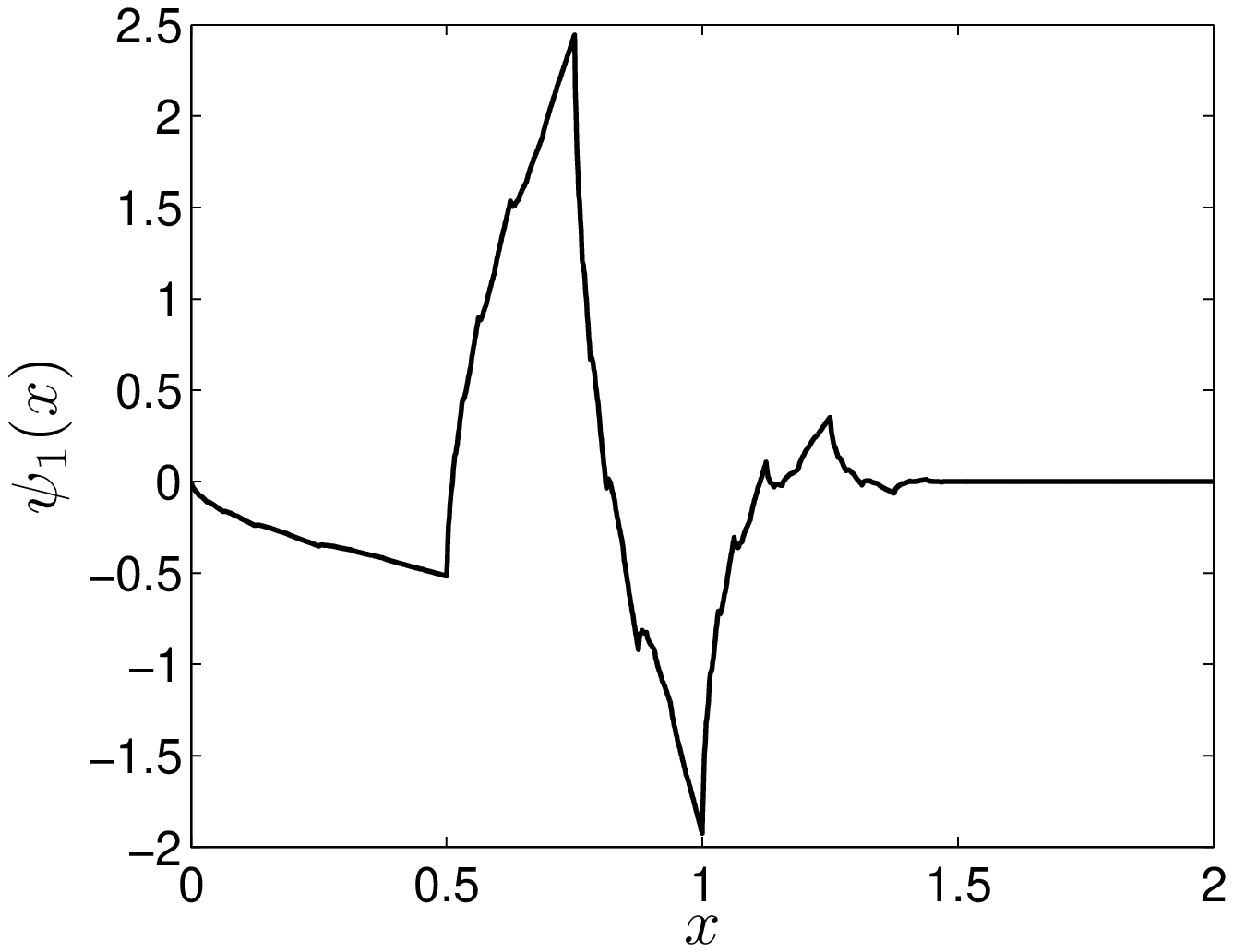}
\par\end{centering}

}\tabularnewline
\subfloat[]{\begin{centering}
\includegraphics[width=2in]{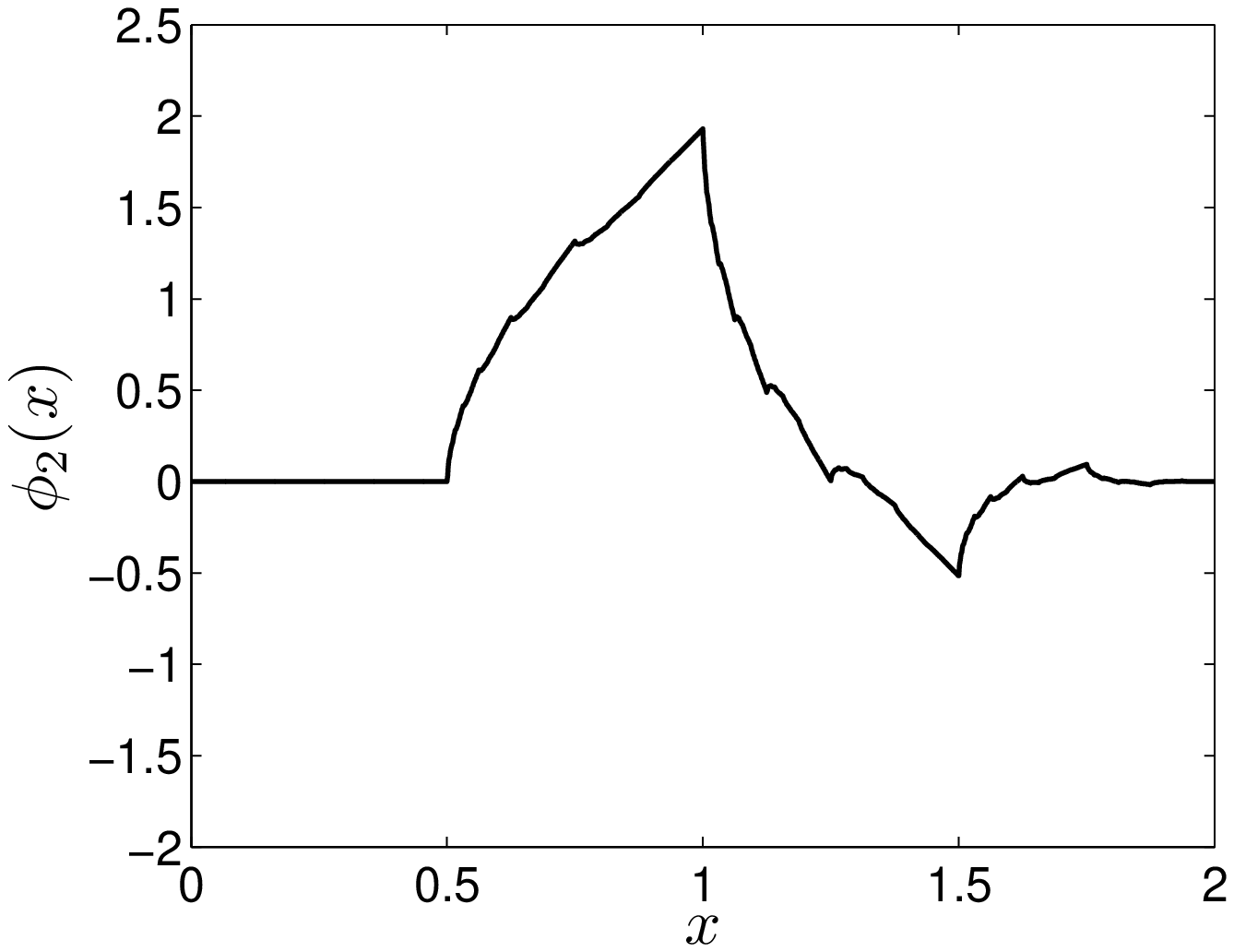}
\par\end{centering}

} & \subfloat[]{\begin{centering}
\includegraphics[width=2in]{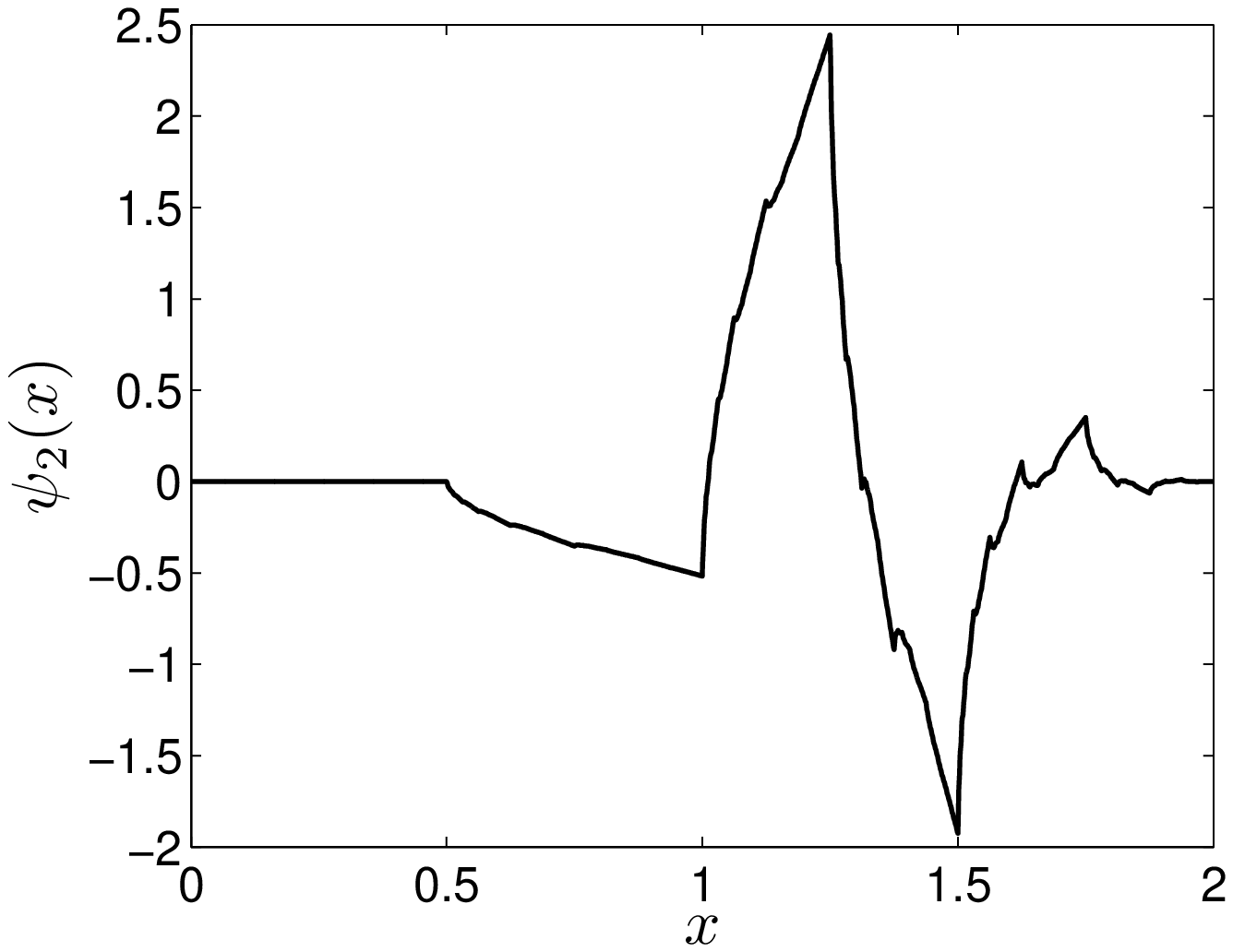}
\par\end{centering}

}\tabularnewline
\end{tabular}
\par\end{centering}

\caption{The standard Daubechies wavelet of order 2 is on the first row to
illustrate that the balanced Daubechies multiwavelet (last two rows)
is just compressed and translated versions of the standard Daubechies
wavelet.\label{fig:BalDaub2}}
\end{figure}

As with wavelets, there are many multiwavelet families---Shen-Tan-Tham
(STT) by \citet{STT2000}, Donovan-Geronimo-Hardin-Massopust (DGHM)
by \citet{GHM1994} and \citet{DGHM1996} (with multiwavelet functions
by \citet{StrangStrela1995}), and Chui-Lan (CL) by \citet{Chui1995},
just to name a few. Figs. \ref{fig:DGHM}, \ref{fig:STT}, and \ref{fig:BalDaub2}
illustrate the multiscaling and multiwavelet pairs for several well-known
multiwavelet families. Like with standard wavelets, the properties
of orthogonality and compact support are enforced during the construction
of the recursion coefficients $H_{k}$ uniquely defining a multiwavelet.
Unlike wavelets, however, multiwavelets can be simultaneously symmetric
and compactly supported while retaining orthogonality across their
integer translates. This desirable property serves as the primary
motivation for density estimation with a multiwavelet basis rather
than the standard wavelet basis.

An interesting incarnation of multiwavelets are the so-called ``balanced
multiwavelets,'' first explicitly introduced by \citet{Lebrun1997}
and later more-rigidly formalized in \citet{Lebrun1998}. \citet{Selesnick1998}
investigated the approximation properties of balanced multiwavelets,
and later \citet{Lebrun2001} generalized the concept of multiwavelet
balancing to arbitrarily high order. This class of multiwavelets was
developed to solve a problem in signal-processing arising from the
vector nature of multiwavelets, where a one-dimensional input signal
must be vectorized before being passed through a multifilter. This
vectorization can result in undesirable effects on the signal reconstruction
due to ``unbalanced'' channels in the lowpass coefficients of the
multifilter. Through procedures proposed by \citet{Lebrun2001}, a
multiwavelet can be balanced to eliminate these undesirable features
of the lowpass coefficients. In fact, standard wavelets, such as the
Daubechies, Symlets, and Coiflets, can be used to construct balanced
multiwavelets. A particular balanced multiwavelet is shown in Fig.
\ref{fig:BalDaub2}. Here, we have taken the standard Daubechies wavelet
of order 2 and, with the toolkit from \citet{Keinert04}, used it
to construct a balanced multiwavelet of multiplicity $r=2$. The resulting
balanced multiwavelet is simply a compressed version of the original
wavelet translated on the half-integers. We have illustrated this
explicitly in Fig. \ref{fig:BalDaub2} by showing the standard Daubechies
wavelet on the first row and the components of the balanced multiwavelet
on the remaining rows. We specifically mention balanced multiwavelets
in this paper as they produce some interesting, though unsurprising,
results when used as a basis for MWDE. 

\section{Multiwavelet Density Estimation\label{sec:MWDE}}

Our objective is to approximate a density function $p(x)$ using the
multiwavelet basis in a form analogous to WDE. Again, the input is
an i.i.d. sample of one-dimensional data $X=\left\{ X_{i}\right\} _{i=1}^{N}$,
and we aim to construct
\begin{equation}
p\left(x\right)=\sum_{k\in\mathbb{Z}}\underline{\alpha}_{j_{0},k}^{T}\underline{\phi}_{j_{0},k}\left(x\right)+\sum_{j=j_{0}}^{\infty}\sum_{k\in\mathbb{Z}}\underline{\beta}_{j,k}^{T}\underline{\psi}_{j,k}\left(x\right),\label{eq:MultiwaveletExpansion}
\end{equation}
where, analogous to the wavelet case, $\underline{\varphi}_{j,k}\left(x\right)=2^{j/2}\underline{\varphi}\left(2^{j}x-k\right)$,
where $\underline{\varphi}$ is $\underline{\phi}$ or $\underline{\psi}$.
The coefficients $\alpha_{j_{0},k}$ and $\beta_{j,k}$ have become
the $r$-dimensional vectors $\underline{\alpha}_{j_{0},k}$ and $\underline{\beta}_{j,k}$.
We can expand Eq. \ref{eq:MultiwaveletExpansion} into its explicit
vector form to see the reconstruction more clearly: 
\begin{equation}
p\left(x\right)=\sum_{k\in\mathbb{Z}}\left(\begin{array}{c}
\alpha_{1,j_{0},k}\\
\vdots\\
\alpha_{r,j_{0},k}
\end{array}\right)^{T}\left(\begin{array}{c}
\phi_{1,j_{0},k}\left(x\right)\\
\vdots\\
\phi_{r,j_{0},k}\left(x\right)
\end{array}\right)+\sum_{j=j_{0}}^{\infty}\sum_{k\in\mathbb{Z}}\left(\begin{array}{c}
\beta_{1,j,k}\\
\vdots\\
\beta_{r,j,k}
\end{array}\right)^{T}\left(\begin{array}{c}
\psi_{1,j,k}\left(x\right)\\
\vdots\\
\psi_{r,j,k}\left(x\right)
\end{array}\right).
\end{equation}
Evidently, the density function is completely described by the coefficients
$\underline{\alpha}_{j_{0},k}$ and $\underline{\beta}_{j,k}$, so
the objective is to estimate $\underline{\alpha}_{j_{0},k}$ and $\underline{\beta}_{j,k}$
as $\hat{\underline{\alpha}}_{j_{0},k}$ and $\hat{\underline{\beta}}_{j,k}$,
respectively, using only the sample $X$. Before we detail a projection
approach similar to WDE, it is worth expanding on the notion of orthonormality
as it applies to multiwavelets to make explicitly clear the idea of
multiwavelet density estimation in an OSE environment. 

Multiwavelets are orthonormal across integer translates $k$ if

\begin{equation}
\left\langle \underline{\phi}\left(x\right),\underline{\phi}\left(x-k\right)\right\rangle =\int\underline{\phi}\left(x\right)\underline{\phi}^{*}\left(x-k\right)dx=\delta_{0k}I,\label{eq:MultiwaveletOrthogonality}
\end{equation}
where $\delta_{ij}$ is the Kronecker delta function, $I$ is the
$r\times r$ identity matrix, and $\underline{\phi}^{*}$ denotes
the conjugate transpose of the vector $\underline{\phi}$. Note, that
since $\underline{\phi}\in\mathbb{R}^{r}$ for our purposes, $\underline{\phi}^{*}=\underline{\phi}^{T}$.
So, expanding Eq. \ref{eq:MultiwaveletOrthogonality}, we find that

\begin{equation}
\int\left(\begin{array}{c}
\phi_{1}\left(x\right)\\
\vdots\\
\phi_{r}\left(x\right)
\end{array}\right)\left(\begin{array}{ccc}
\phi_{1}\left(x-k\right) & \cdots & \phi_{r}\left(x-k\right)\end{array}\right)dx=\delta_{0k}\left(\begin{array}{ccc}
1 &  & 0\\
 & \ddots\\
0 &  & 1
\end{array}\right),
\end{equation}
which implies
\begin{equation}
\left(\begin{array}{ccc}
\int\phi_{1}\left(x\right)\phi_{1}\left(x-k\right)dx & \cdots & \int\phi_{1}\left(x\right)\phi_{r}\left(x-k\right)dx\\
\vdots & \ddots & \vdots\\
\int\phi_{r}\left(x\right)\phi_{1}\left(x-k\right)dx & \cdots & \int\phi_{r}\left(x\right)\phi_{r}\left(x-k\right)dx
\end{array}\right)=\delta_{0k}\left(\begin{array}{ccc}
1 &  & 0\\
 & \ddots\\
0 &  & 1
\end{array}\right).
\end{equation}
Finally, it is evident that
\begin{equation}
\int\phi_{i}\left(x\right)\phi_{j}\left(x-k\right)dx=\delta_{0k}\delta_{ij},
\end{equation}
implying $\phi_{i}\left(x\right)$ is orthonormal to $\phi_{j}\left(x\right)$
when $i\neq j$ across integer translates $k$; that is, $\underline{\phi}$
can be constructed such that its elements are orthonormal across integer
translates, justifying Eq. \ref{eq:MultiwaveletExpansion}. The same
conclusion holds for the mother multiwavelet functions. Therefore,
we can calculate the coefficients $\underline{\alpha}_{j_{0},k}$
and $\underline{\beta}_{j,k}$ using the standard inner product projection:
\begin{equation}
\underline{\alpha}_{j_{0},k}=\left\langle p,\underline{\phi}_{j_{0},k}\right\rangle =\int p\left(x\right)\underline{\phi}_{j_{0},k}\left(x\right)dx,\label{eq:MWinnerProduct}
\end{equation}
and similarly for $\underline{\beta}_{j,k}$:
\[
\underline{\beta}_{j,k}=\left\langle p,\underline{\psi}_{j,k}\right\rangle =\int p\left(x\right)\underline{\psi}_{j,k}\left(x\right)dx.
\]
As before, this allows us to interpret the inner product in Eq. \ref{eq:MWinnerProduct}
as an expectation
\begin{equation}
\underline{\alpha}_{j_{0},k}=\int p\left(x\right)\underline{\phi}_{j_{0},k}\left(x\right)dx=\mathcal{E}\left[\underline{\phi}_{j_{0},k}\left(x\right)\right],
\end{equation}
which is approximated as the sample mean
\begin{equation}
\hat{\underline{\alpha}}_{j_{0},k}=\frac{1}{N}\sum_{i=1}^{N}\underline{\phi}_{j_{0},k}\left(X_{i}\right).
\end{equation}
The multiwavelet function coefficients $\underline{\beta}_{j,k}$
are estimated as expected:
\begin{equation}
\hat{\underline{\beta}}_{j,k}=\frac{1}{N}\sum_{i=1}^{N}\underline{\psi}_{j,k}\left(X_{i}\right).
\end{equation}
The final approximation $\hat{p}\left(x\right)$ is given by
\begin{equation}
\hat{p}\left(x\right)=\sum_{k\in\mathbb{Z}}\hat{\underline{\alpha}}_{j_{0},k}^{T}\underline{\phi}_{j_{0},k}\left(x\right)+\sum_{j=j_{0}}^{J}\sum_{k\in\mathbb{Z}}\underline{\hat{\beta}}_{j,k}^{T}\underline{\psi}_{j,k}\left(x\right).
\end{equation}
This is the linear multiwavelet density estimator. As in the wavelet
case, nonlinear MWDE is also possible using either hard or soft thresholding
of the multiwavelet coefficients $\hat{\underline{\beta}}_{j,k}$.
For instance, see the work already done by \citep{Downie1998,StrelaAndWalden1998,Strela1999}
on multiwavelet coefficient thresholding in signals-processing applications.
In addition to element-wise thresholding, \citet{Bacchelli2002} investigated
vector-wise thresholding, taking advantage of the fact that multiwavelet
coefficients can be correlated in their vector representations. As
in WDE, the resolution levels $j_{0}\leq j\leq J<\infty$ are chosen
by the user or by using model selection techniques as described, for
instance, by \citet{Vidakovic1999}. We do not rigorously address
the choice of resolution levels or coefficient thresholding here,
as our objective is only to show that OSE can be successfully accomplished
with multiwavelet bases.

\section{Experimental Results\label{sec:Experiments}}

We detail several experiments using MWDE in comparison with standard
WDE. We examine how multiwavelet symmetry can affect the reconstruction
of densities with global and/or local symmetrical properties. Finally,
we investigate the interesting case of using balanced multiwavelets.
In addition to simple densities like the Gaussian, we extensively
use the complicated density functions covered in \citet{Marron92}
and \citet{Wand93}. These functions are, in many cases, multi-modal
and contain local symmetries. We constructed a density estimator using
various multiwavelet bases and measured the accuracy of the estimator
under the integrated square error ($\mathrm{ISE}$) between the estimated
density $\hat{p}$ and the actual density $p$. For the illustrated
test cases, we use a variety of multiwavelet families chosen by their
regularity of appearance in the literature. The multiwavelets were
computationally constructed using the iterative cascading algorithm
\citep{Strang97}. This is a standard procedure used even for wavelet
constructions, with implementation details specific to multiwavelets
available in \citet{Keinert04}.

\subsection{Multiwavelet Symmetry}

We begin by attempting to empirically motivate the advantages of MWDE
versus WDE. Theoretically this was based on the fact that multiwavelets
can be constructed with simultaneous symmetry, compact support, and
orthogonality, whereas wavelets can not have these properties simultaneously.
We developed a simple test to investigate the utility of the symmetry
property by using multiwavelets and standard wavelets to estimate
the symmetric unimodal Gaussian and a bimodal distribution. For this
comparison, we employed the commonly-used Daubechies wavelet for the
WDE. For multiwavelets, we have a choice of several symmetric and
symmetric/antisymmetric families that have already been developed.
The DGHM multiwavelet \citep{DGHM1996} has symmetric multiscaling
functions and symmetric/antisymmetric multiwavelet functions \citep{StrangStrela1995};
it is plotted in Fig. \ref{fig:DGHM}. The STT multiwavelet \citep{STT2000}
has symmetric/antisymmetric multiscaling and multiwavelet functions
and is plotted in Fig. \ref{fig:STT}. We use both the DGHM and STT
multiwavelets to demonstrate the capabilities of symmetric multiwavelets
for density estimation. We do not employ any MRA for our density estimation,
so the multiwavelet functions are not of direct interest here. However,
such a density estimator could be easily constructed following the
methodology detailed in \S\ref{sec:MWDE}.

\begin{figure}
\noindent \begin{centering}
\begin{tabular}{cc}
\subfloat[$j=-2$]{\begin{centering}
\includegraphics[width=2in]{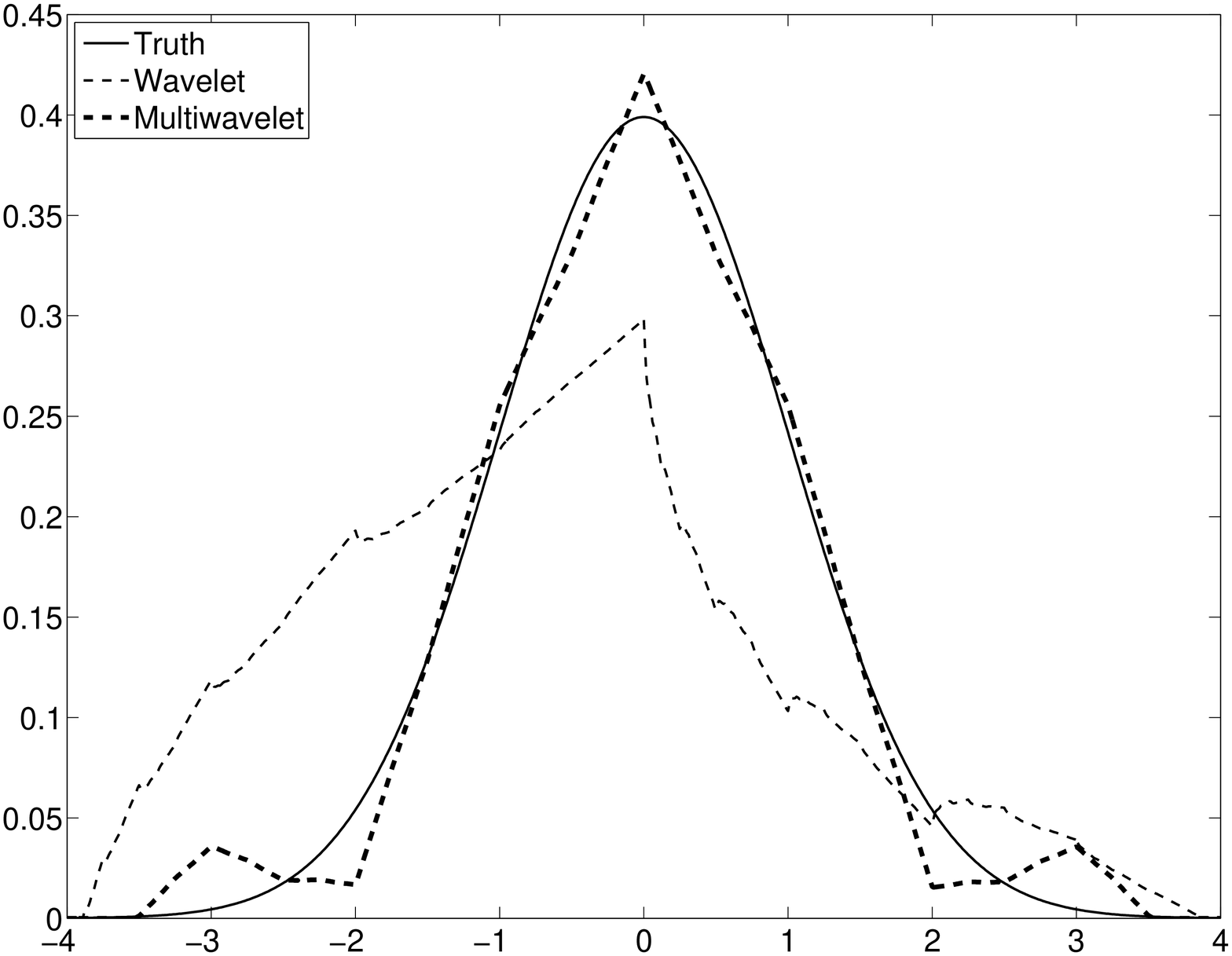}
\par\end{centering}

} & \subfloat[$j=-1$]{\begin{centering}
\includegraphics[width=2in]{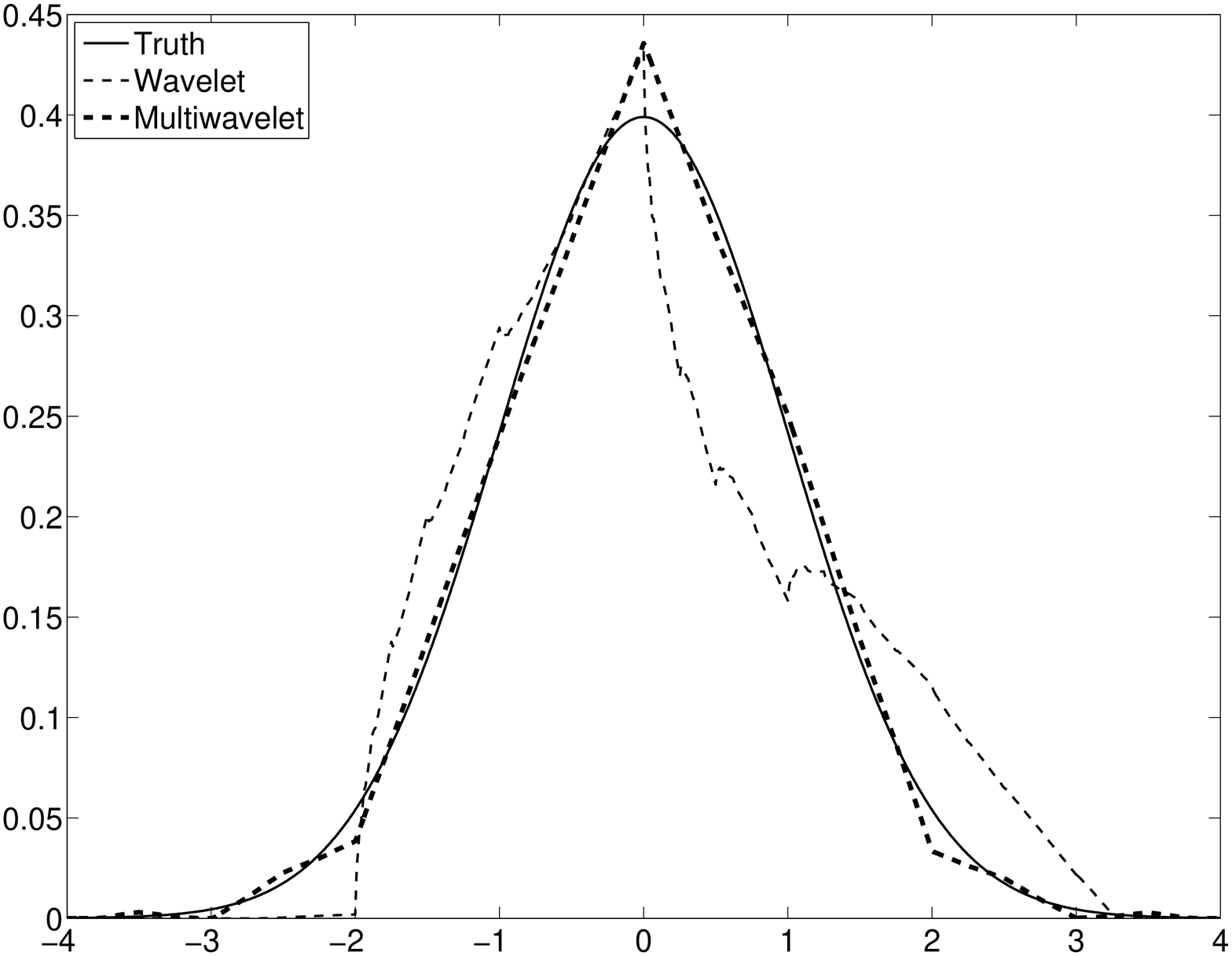}
\par\end{centering}

}\tabularnewline
\subfloat[$j=0$]{\noindent \centering{}\includegraphics[width=2in]{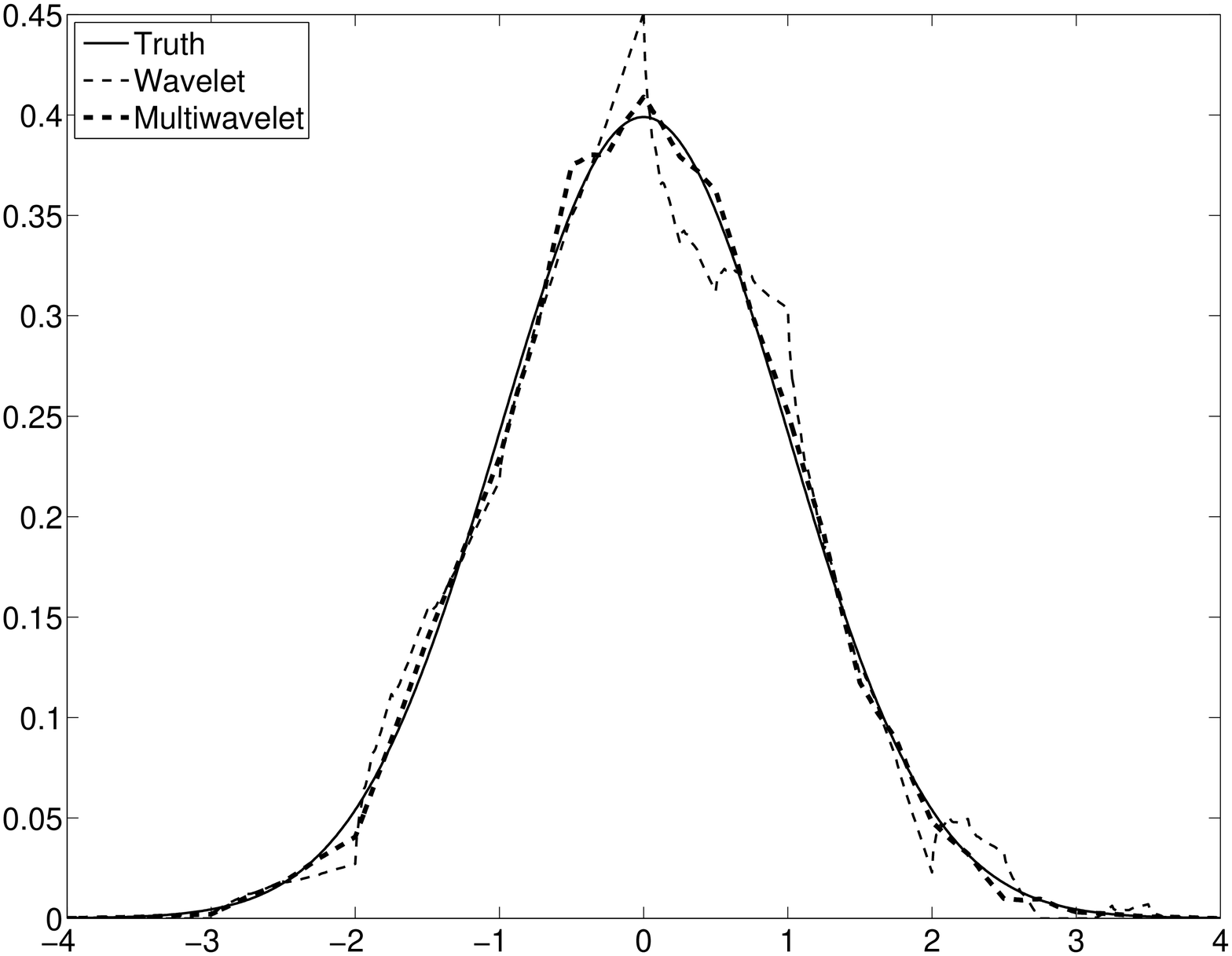}} & \subfloat[$j=1$]{\noindent \begin{centering}
\includegraphics[width=2in]{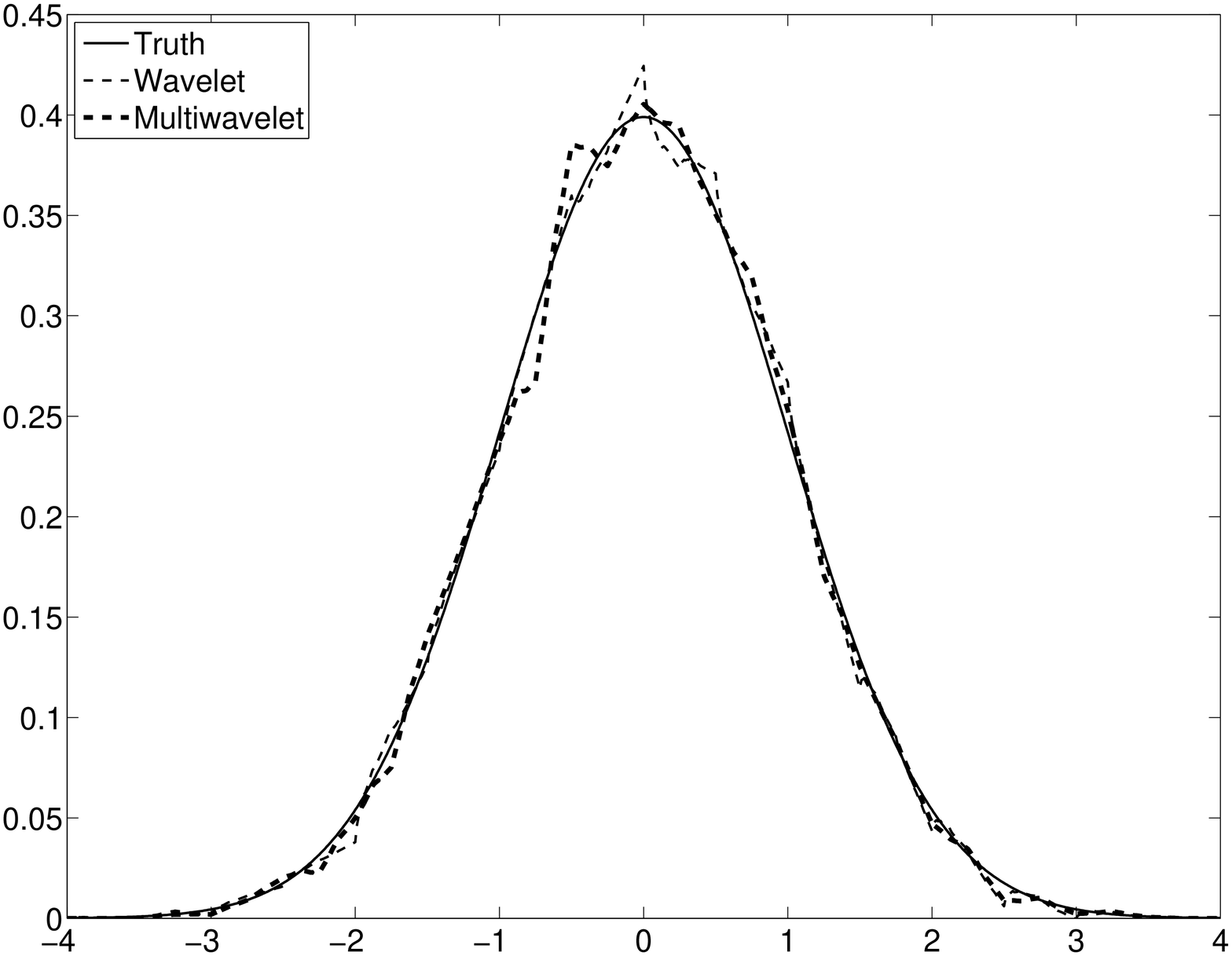}
\par\end{centering}

}\tabularnewline
\subfloat[$j=2$]{\begin{centering}
\includegraphics[width=2in]{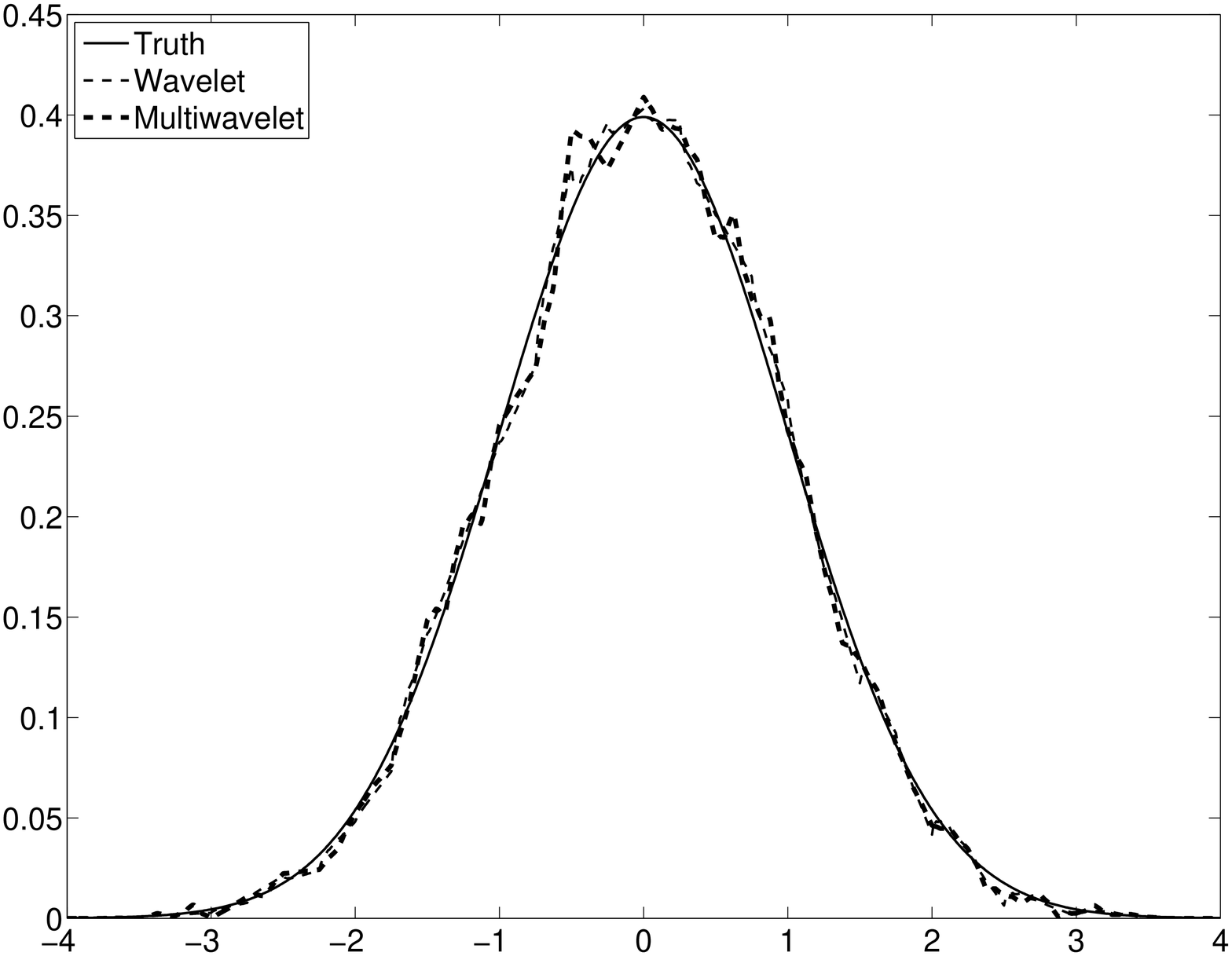}
\par\end{centering}

} & \subfloat[$j=3$]{\begin{centering}
\includegraphics[width=2in]{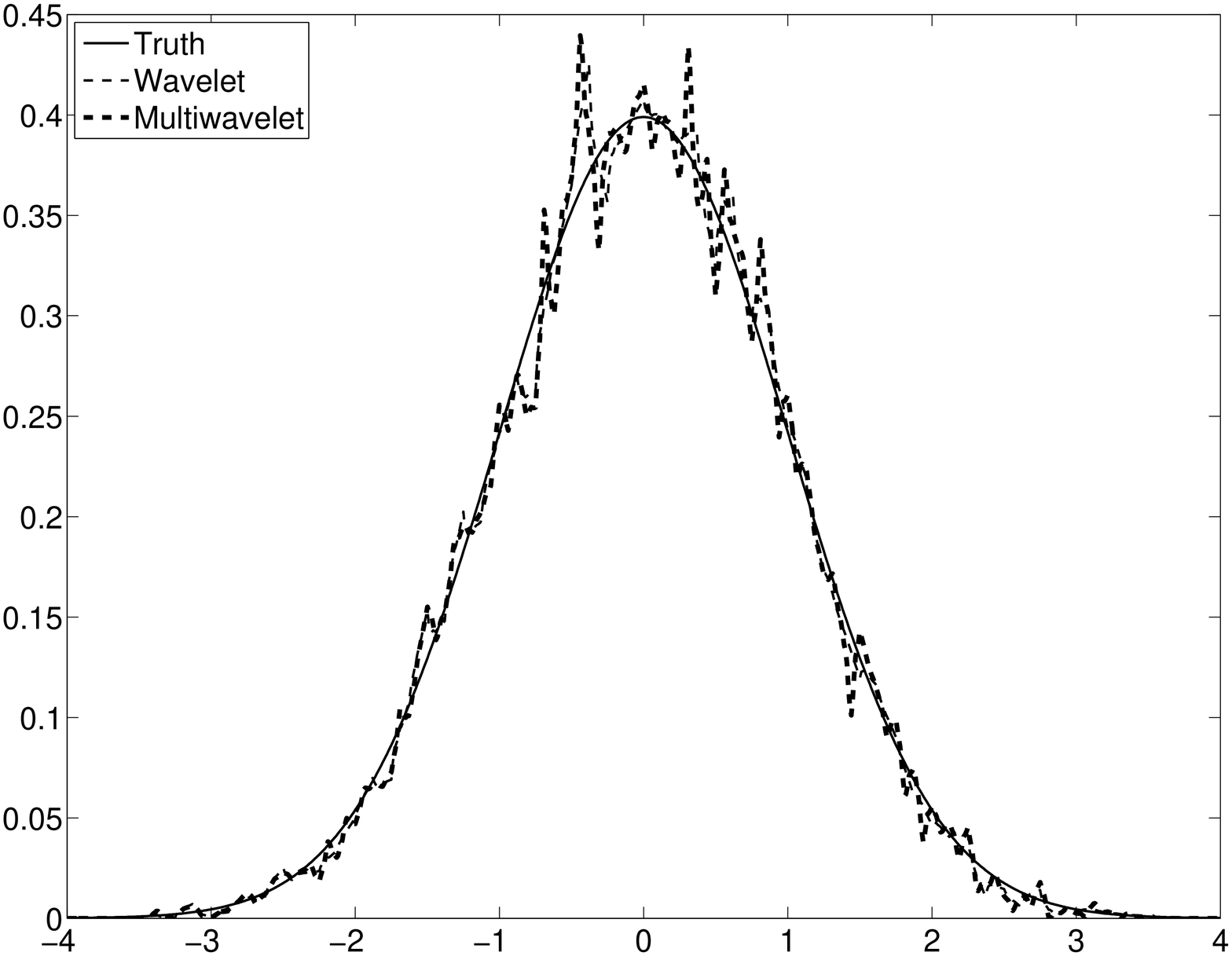}
\par\end{centering}

}\tabularnewline
\end{tabular}
\par\end{centering}

\caption{MWDE of a symmetrical density (Gaussian with mean 0 and standard deviation
1) with 10000 samples across resolutions $-2\leq j\leq3$ with the
symmetric DGHM multiwavelet compared with WDE using the asymmetric
Daubechies wavelet of order 2.\label{fig:DGHMandGauss}}
\end{figure}

First, we estimate a standard Gaussian density, as it has well-known
properties, namely co-located mean/mode/median and symmetry. In Fig.
\ref{fig:DGHMandGauss}, we show the DGHM multiwavelet outperforms
the Daubechies wavelet of order 2 for resolution levels $-2\leq j\leq0$,
but at $j=0$, the Daubechies wavelet begins to produce a comparable
estimate of the density, which continues to improve at finer resolutions.
Given the highly asymmetrical properties of this lower-order Daubechies
wavelet, the MWDE with a basis family such as DGHM was, as expected,
able to outperform the WDE at some coarser resolution levels.

\begin{figure}
\begin{centering}
\begin{tabular}{cc}
\subfloat[$j=-2$]{\begin{centering}
\includegraphics[width=2in]{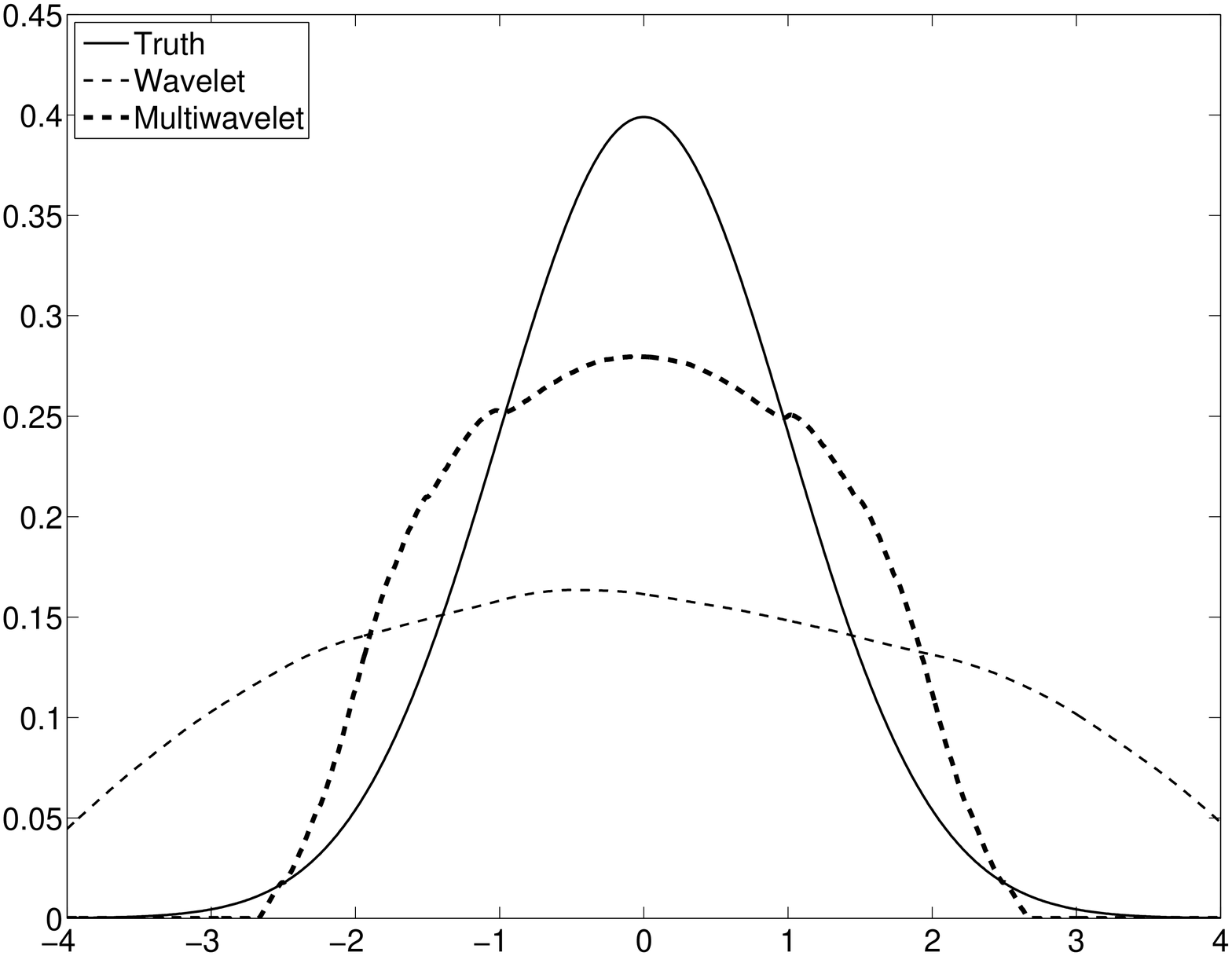}
\par\end{centering}

} & \subfloat[$j=-1$]{\begin{centering}
\includegraphics[width=2in]{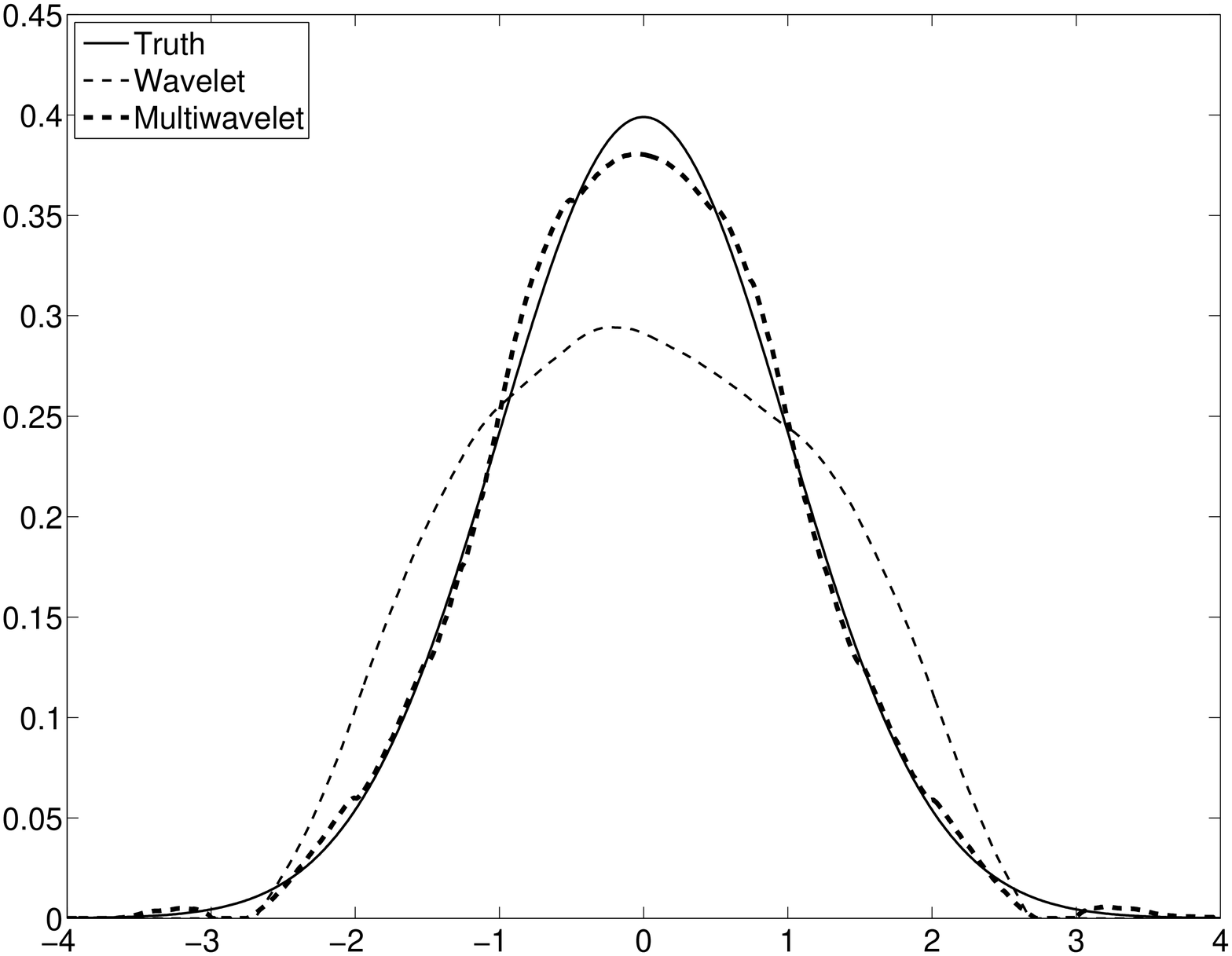}
\par\end{centering}

}\tabularnewline
\subfloat[$j=0$]{\centering{}\includegraphics[width=2in]{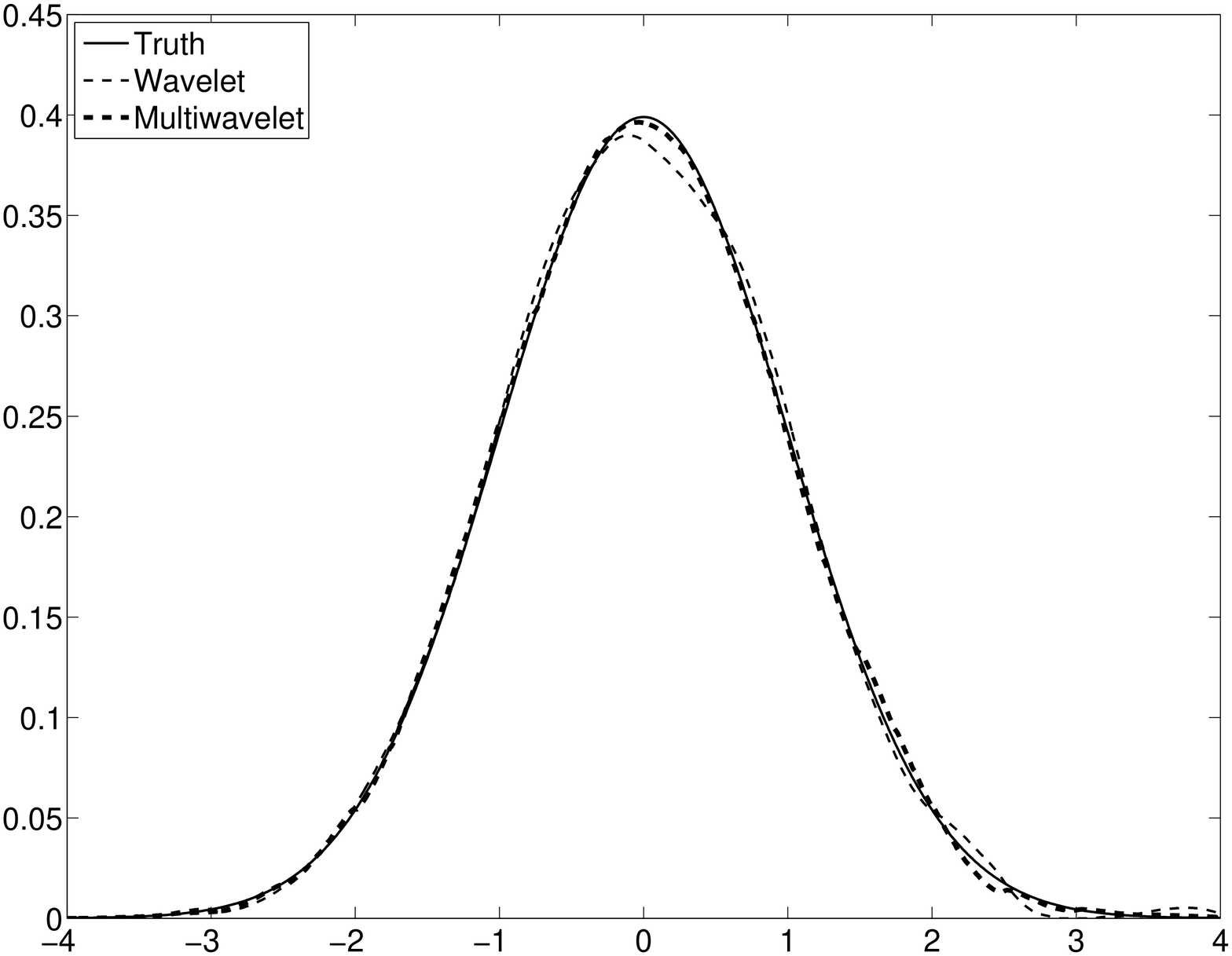}} & \subfloat[$j=1$]{\begin{centering}
\includegraphics[width=2in]{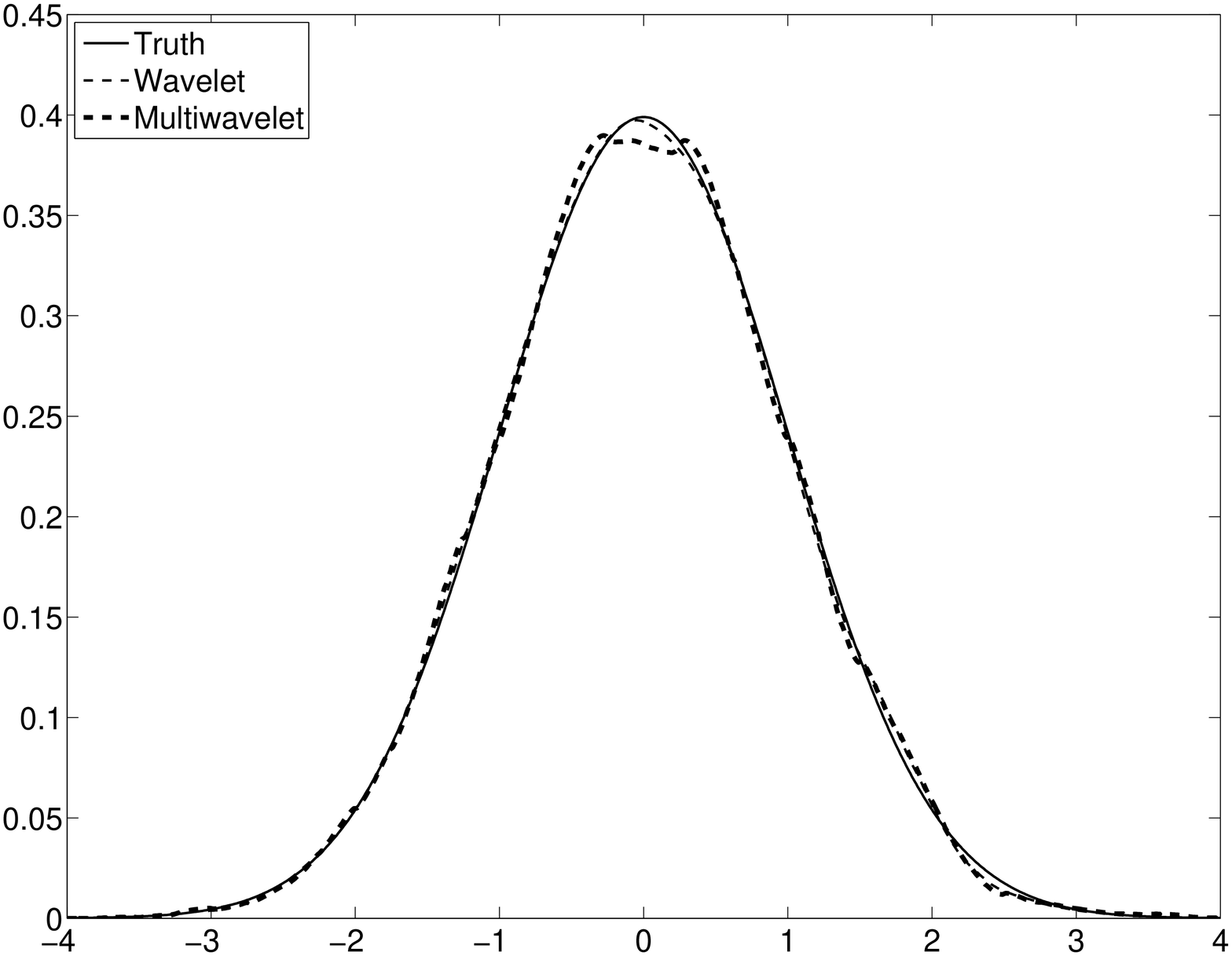}
\par\end{centering}

}\tabularnewline
\subfloat[$j=2$]{\begin{centering}
\includegraphics[width=2in]{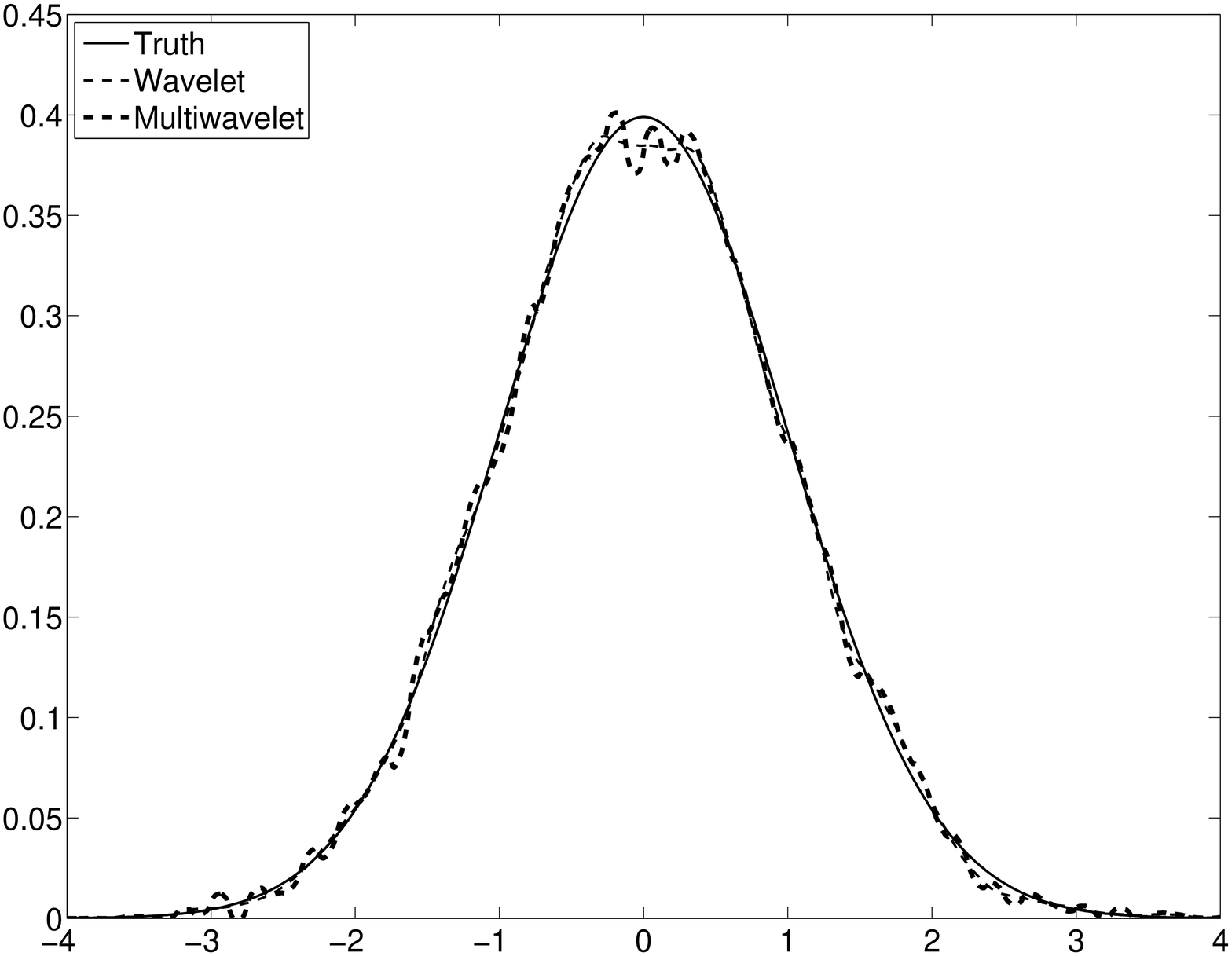}
\par\end{centering}

} & \subfloat[$j=3$]{\begin{centering}
\includegraphics[width=2in]{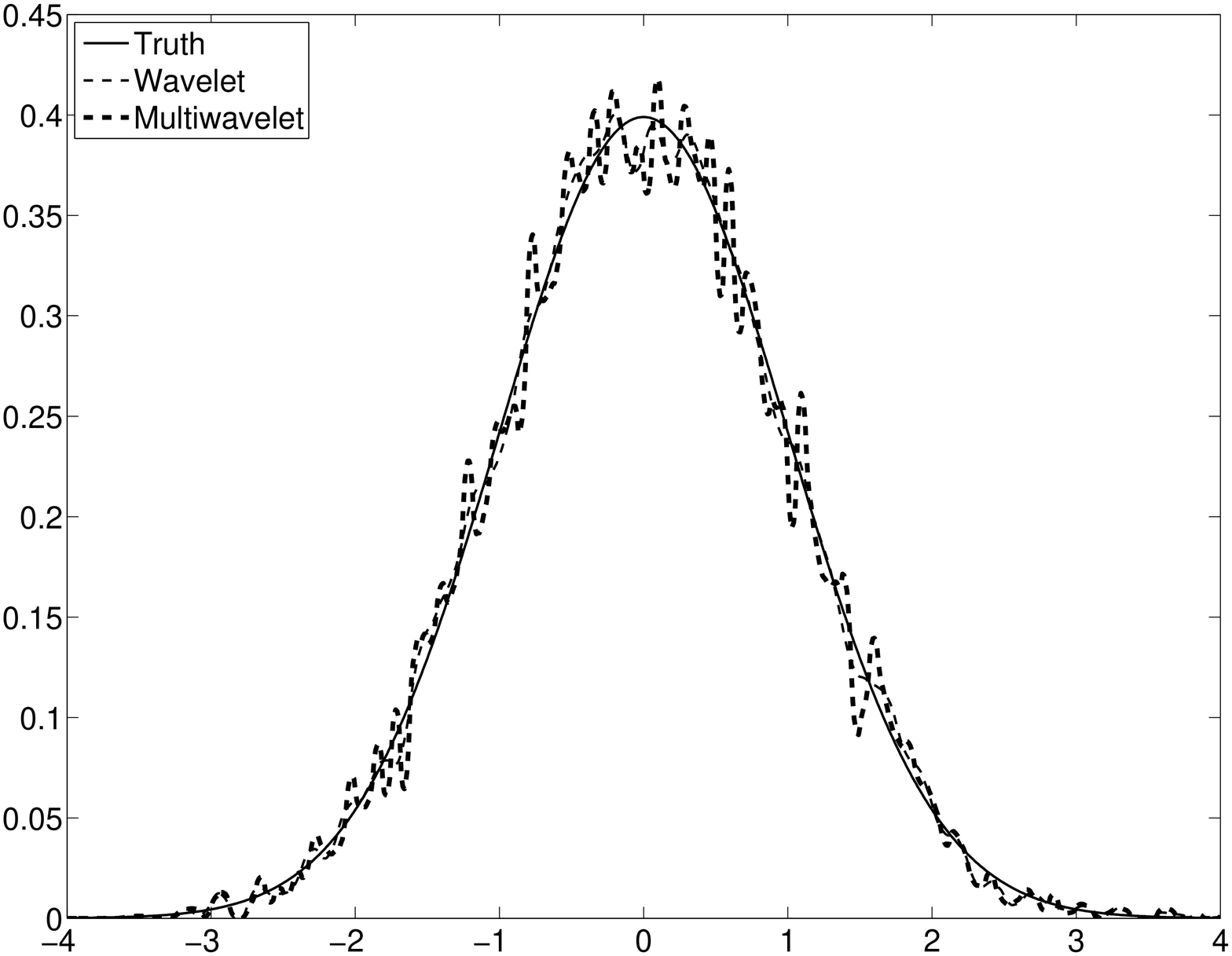}
\par\end{centering}

}\tabularnewline
\end{tabular}
\par\end{centering}

\caption{MWDE of a symmetrical density (Gaussian with mean 0 and standard deviation
1) with 10000 samples across resolutions $-2\leq j\leq3$ with the
symmetric/antisymmetric STT multiwavelet compared with WDE using the
asymmetric Daubechies wavelet of order 2.\label{fig:STTandGauss}}
\end{figure}

It is well-known that the Daubechies wavelet of order 2 is not necessarily
well-suited for estimating smooth densities, such as the Gaussian
distribution; this is evident from the jagged appearance of the wavelet.
Similarly, the DGHM multiwavelet, though symmetrical, is also somewhat
jagged in appearance. A more natural choice for estimating a smooth
and symmetric density would be a higher-order wavelet, such as the
Daubechies wavelet of order 5, and a smoother multiwavelet, such at
the STT multiwavelet. With this, we compare WDE and MWDE using the
more suitable bases just mentioned. The results are shown in Fig.
\ref{fig:STTandGauss}. As anticipated, both the MWDE and the WDE
are superior to the ones in Fig. \ref{fig:DGHMandGauss}. Even when
using the Daubechies wavelet of order 5, we see that MWDE outperforms
WDE at the coarse resolution levels $-2\leq j\leq0$. At finer resolution
levels $j>0$, we see, as in Fig. \ref{fig:DGHMandGauss}, the standard
wavelet basis was able to produce a good density approximation. Being
inherently symmetric, the multiwavelets are able to better reconstruct
the symmetric peaks of these distributions, even at coarse resolutions.
The standard wavelet resolutions must be increased to comparably model
the symmetries. 

\subsection{Balanced Multiwavelets\label{subsec:Validation}}

From the work of \citet{Lebrun1998} and \citet{Lebrun2001}, we know
any wavelet can be used to construct a related balanced multiwavelet
of arbitrarily high multiplicity. For instance, a Daubechies wavelet
can be balanced to a Daubechies multiwavelet of multiplicity $r$,
with $r$ multiscaling and $r$ multiwavelet functions. Balanced multiwavelets
are of particular interest to us because they provide a somewhat-less
subjective method of comparing WDE to MWDE. That is, we can compare
MWDE using a balanced Daubechies multiwavelet to WDE using the Daubechies
wavelet used to produce the balanced multiwavelet. The multiscaling
and multiwavelet functions of balanced Daubechies multiwavelets turn
out to be compressed and translated versions of the corresponding
scaling and wavelet functions of the Daubechies wavelet \citep{Lebrun1998};
this is evident from viewing Fig. \ref{fig:BalDaub2}. 

\begin{figure}
\noindent \begin{centering}
\begin{tabular}{cc}
\subfloat[$j=-2$]{\noindent \begin{centering}
\includegraphics[width=2in]{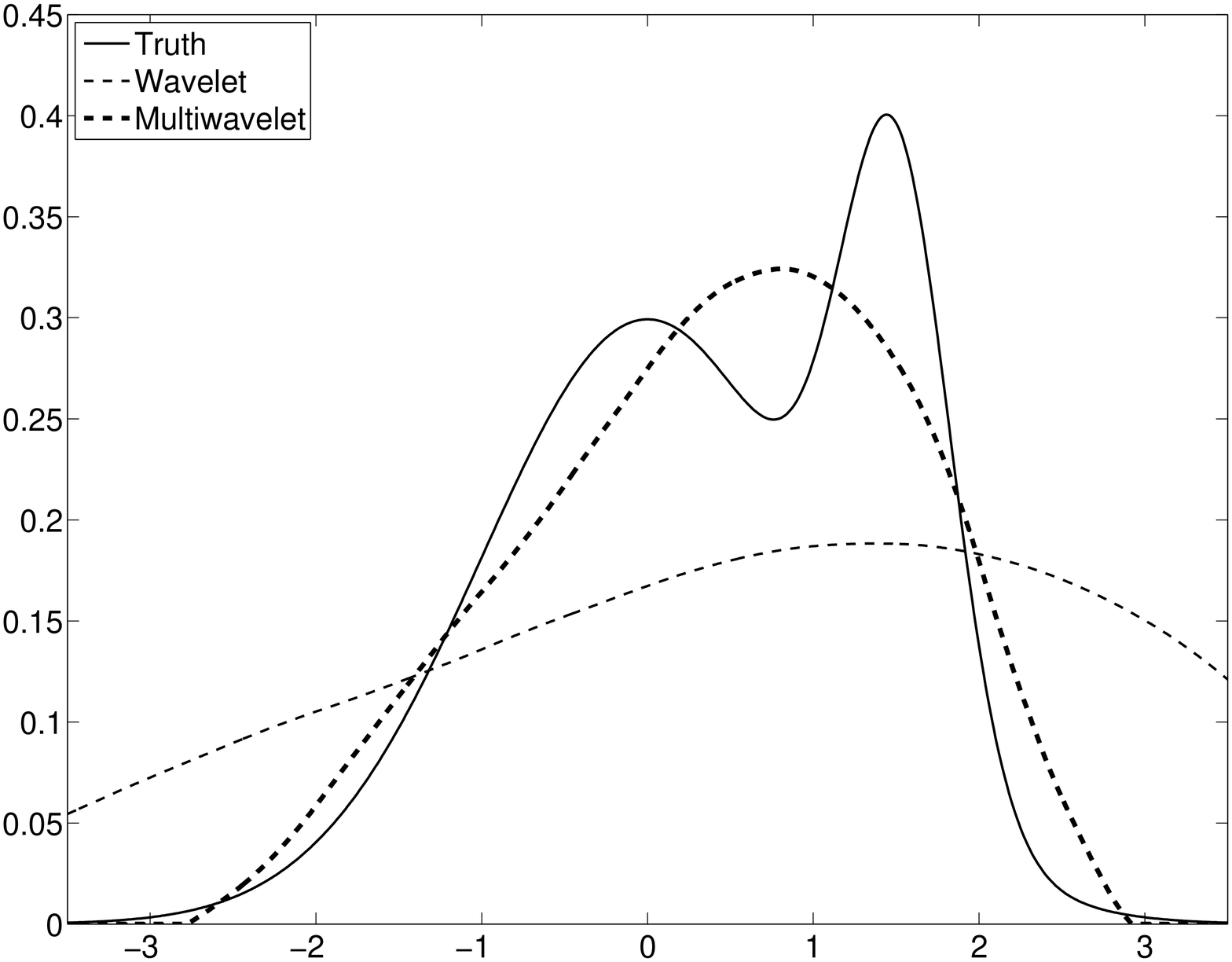}
\par\end{centering}

} & \subfloat[$j=-1$]{\noindent \begin{centering}
\includegraphics[width=2in]{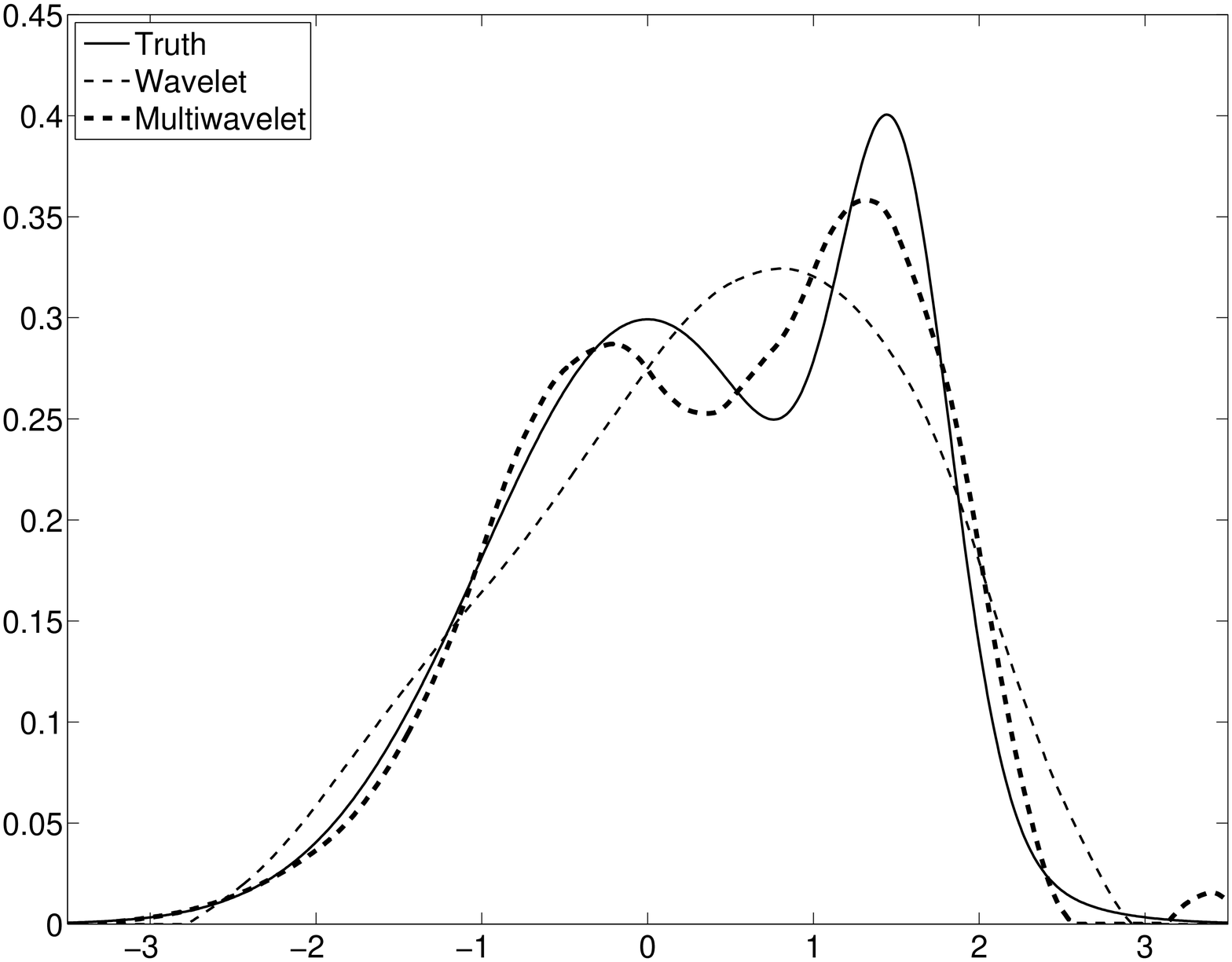}
\par\end{centering}

}\tabularnewline
\subfloat[$j=0$]{\noindent \centering{}\includegraphics[width=2in]{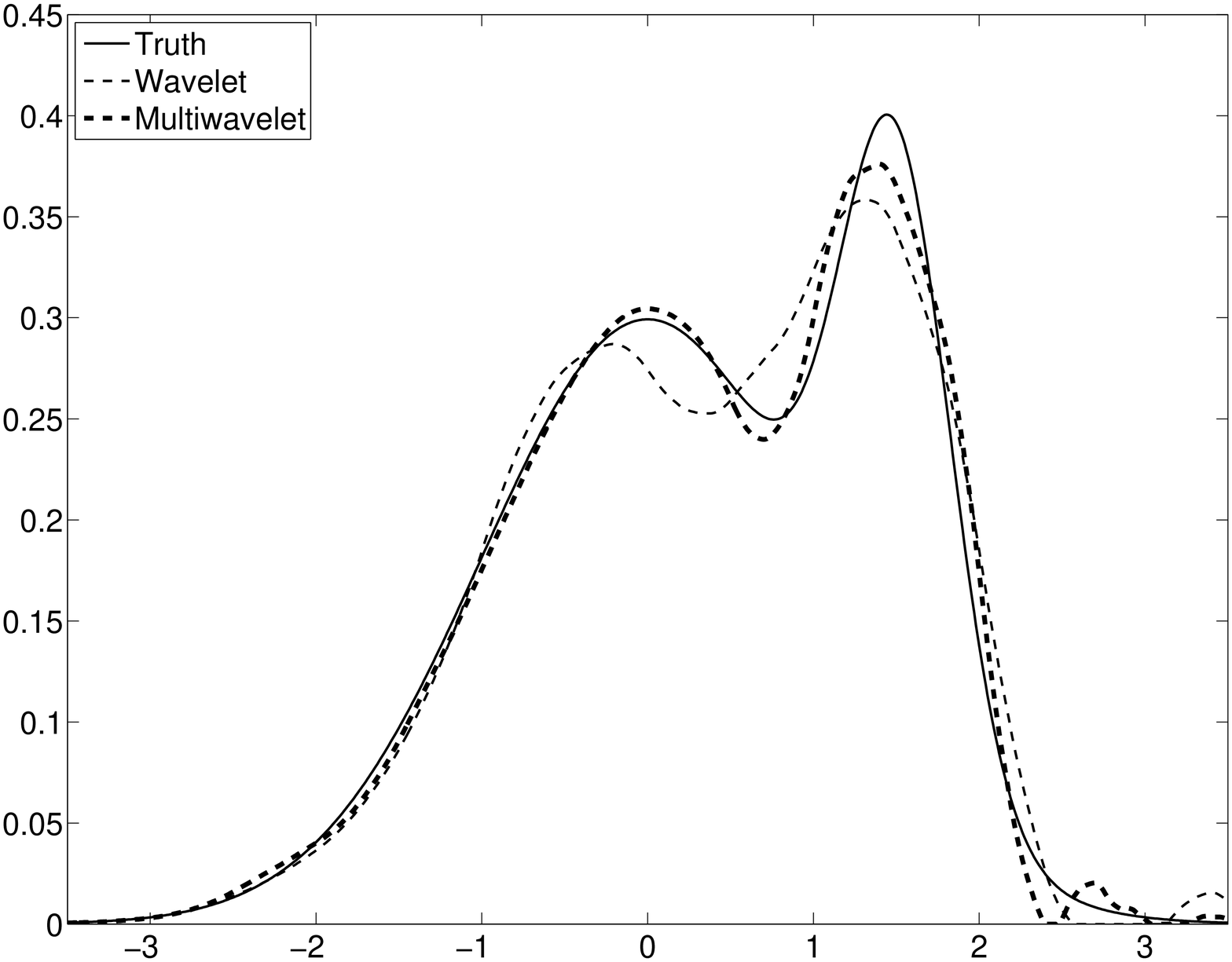}} & \subfloat[$j=1$]{\noindent \begin{centering}
\includegraphics[width=2in]{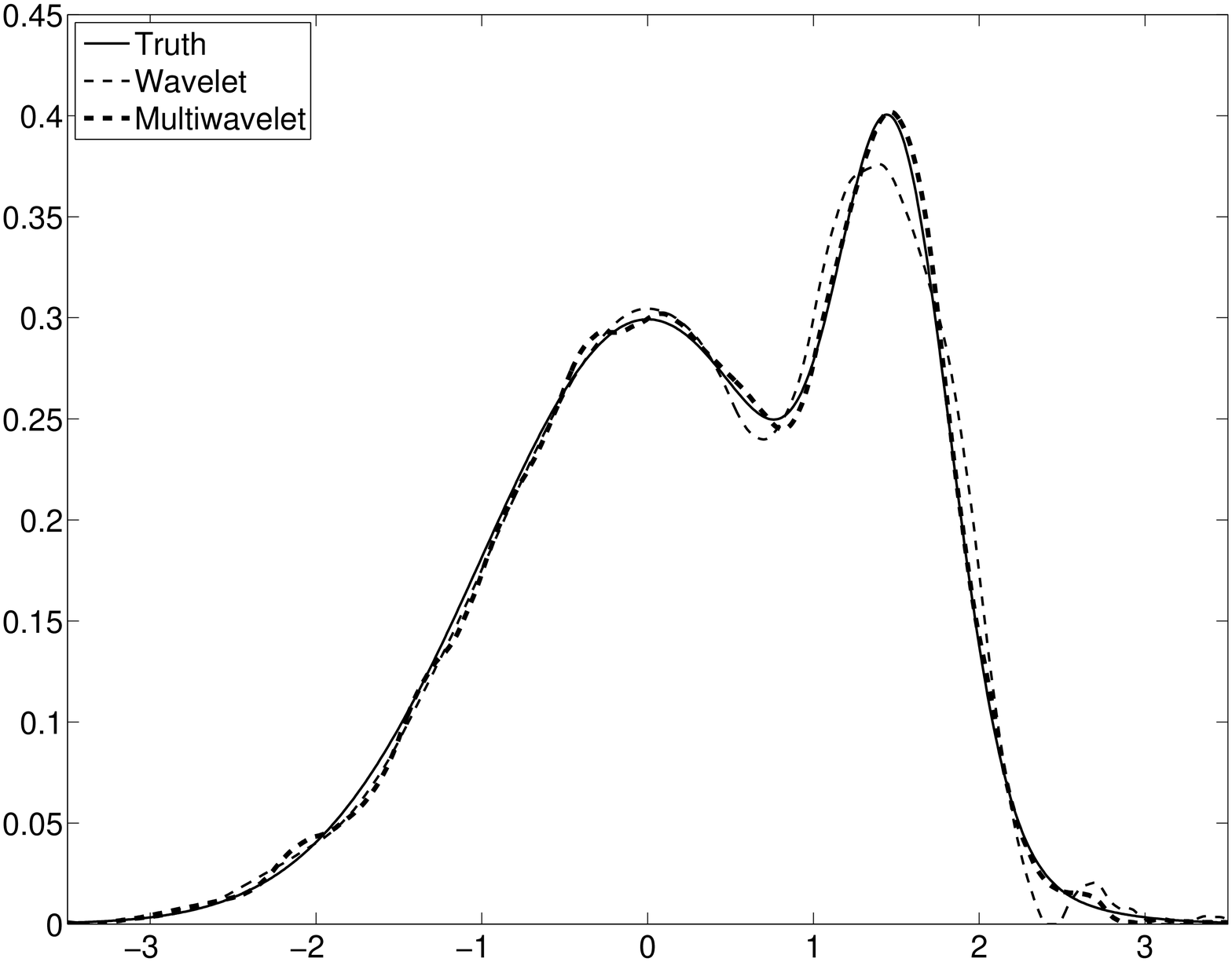}
\par\end{centering}

}\tabularnewline
\subfloat[$j=2$]{\noindent \begin{centering}
\includegraphics[width=2in]{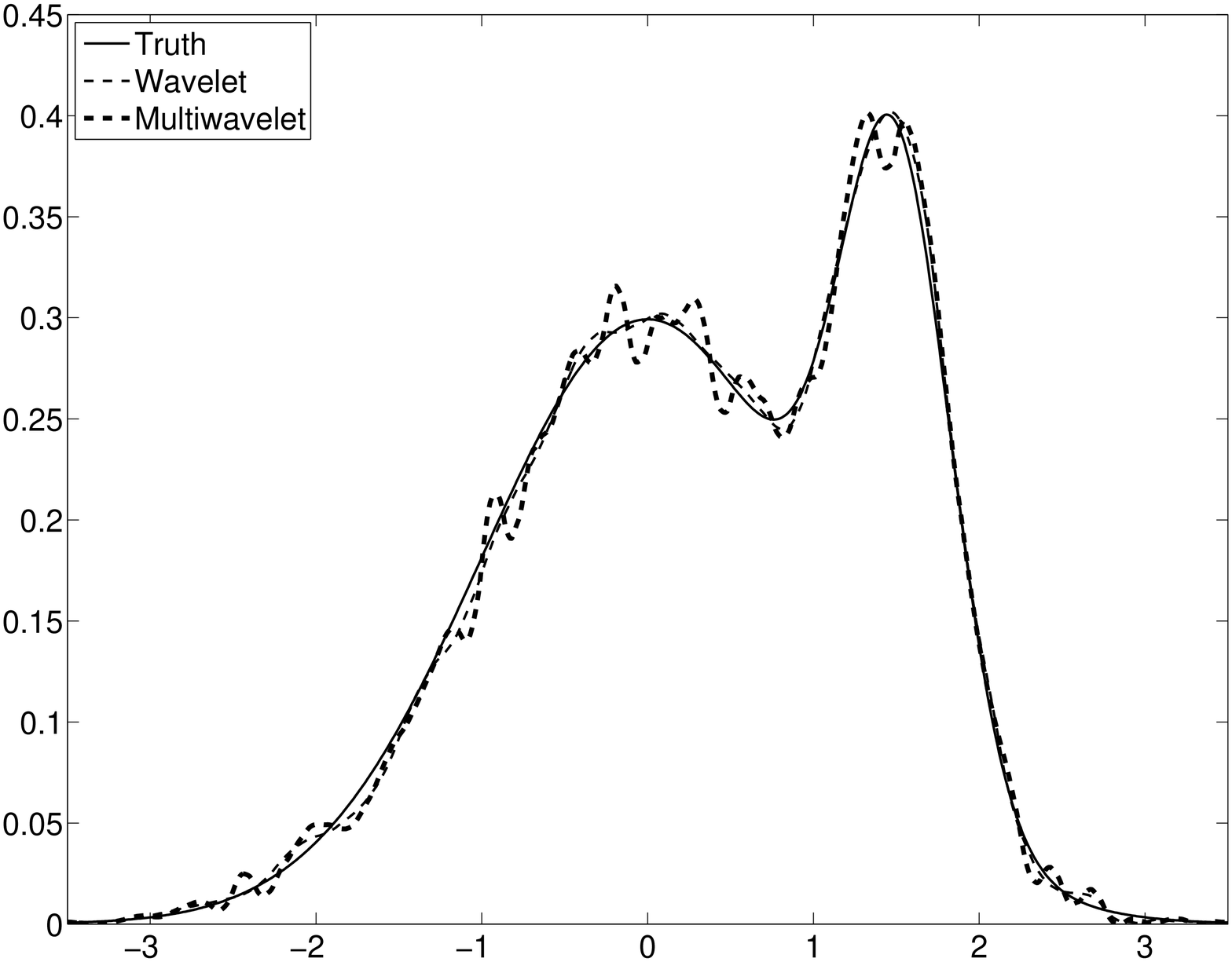}
\par\end{centering}

} & \subfloat[$j=3$]{\noindent \begin{centering}
\includegraphics[width=2in]{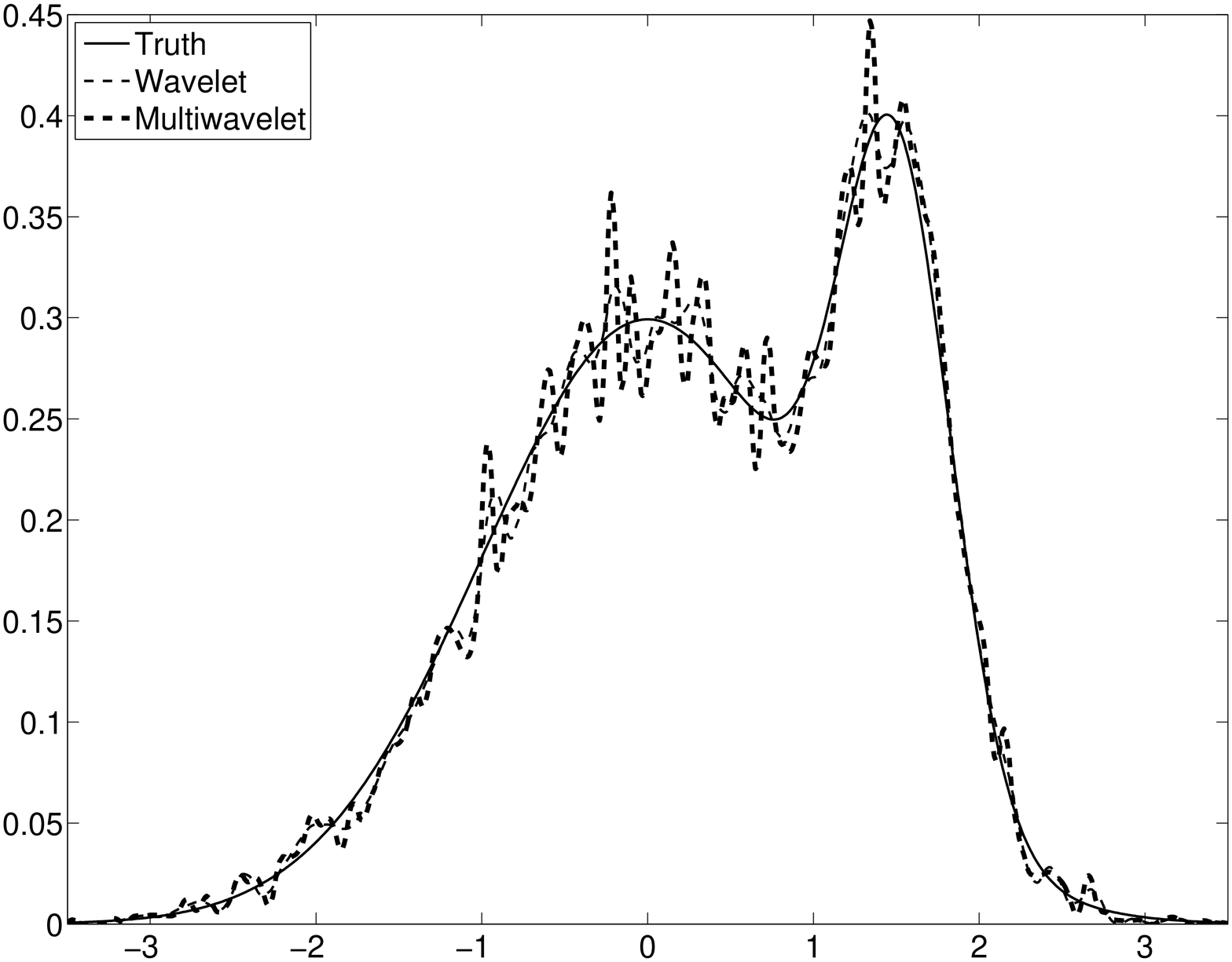}
\par\end{centering}

}\tabularnewline
\end{tabular}
\par\end{centering}

\centering{}\caption{WDE using Daubechies wavelet of order 5 and MWDE using balanced Daubechies
multiwavelet of order 5 compared on a skewed bimodal distribution
with 10000 samples at several resolution levels. Notice how the MWDE
reconstruction at resolution $j=-2$ is very closely comparable to
the WDE reconstruction at $j=-1$. \label{fig:BalProgression}}
\end{figure}

MWDE with balanced multiwavelets possesses some interesting properties,
especially for balanced Daubechies multiwavelets. Because balanced
Daubechies multiwavelets are simply compressed and translated versions
of their wavelet counterparts, we expect MWDE and WDE with these bases
to be closely comparable. In fact, we find exactly this. We empirically
observe a very interesting---if not expected---property where, if
MWDE with balanced multiwavelets of multiplicity $r=2$ produces a
certain reconstruction at some resolution level $j$, then the corresponding
WDE ``lags'' the MWDE, and produces the same, or very close to the
same, reconstruction at resolution $j+1$. This is clearly illustrated
in Fig. \ref{fig:BalProgression}, where we estimate a bimodal distribution
with 10000 samples using the Daubechies wavelet of order 5 for WDE
and the balanced Daubechies multiwavelet of order 5 for MWDE.

\subsection{Other Multiwavelet Families}

\begin{figure}
\begin{centering}
\begin{tabular}{cc}
\subfloat[Balanced Daubechies 6, Daubechies 6, $j=3$]{\begin{centering}
\includegraphics[width=2in]{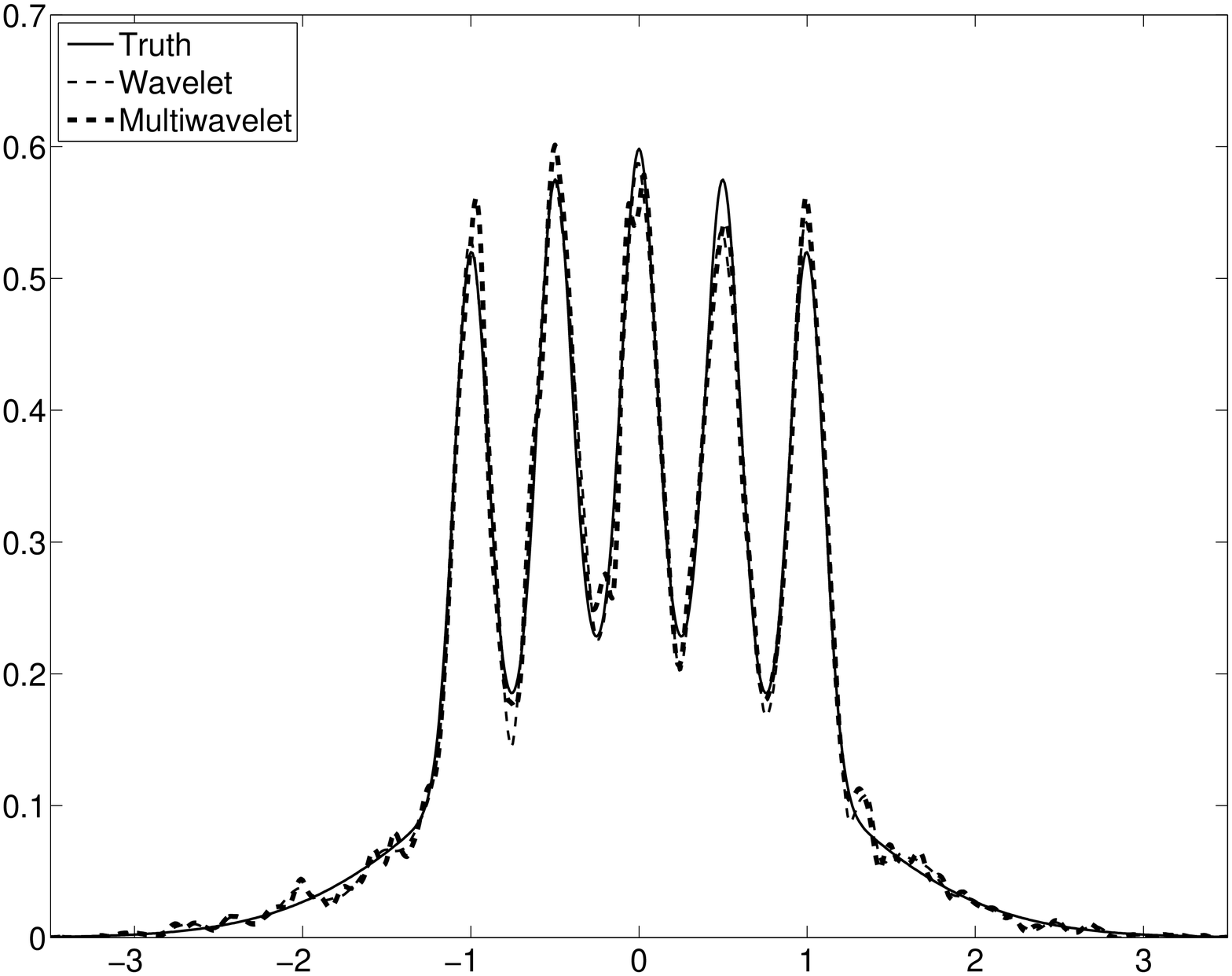}
\par\end{centering}

} & \subfloat[DGHM, Coiflet 5, $j=4$]{\begin{centering}
\includegraphics[width=2in]{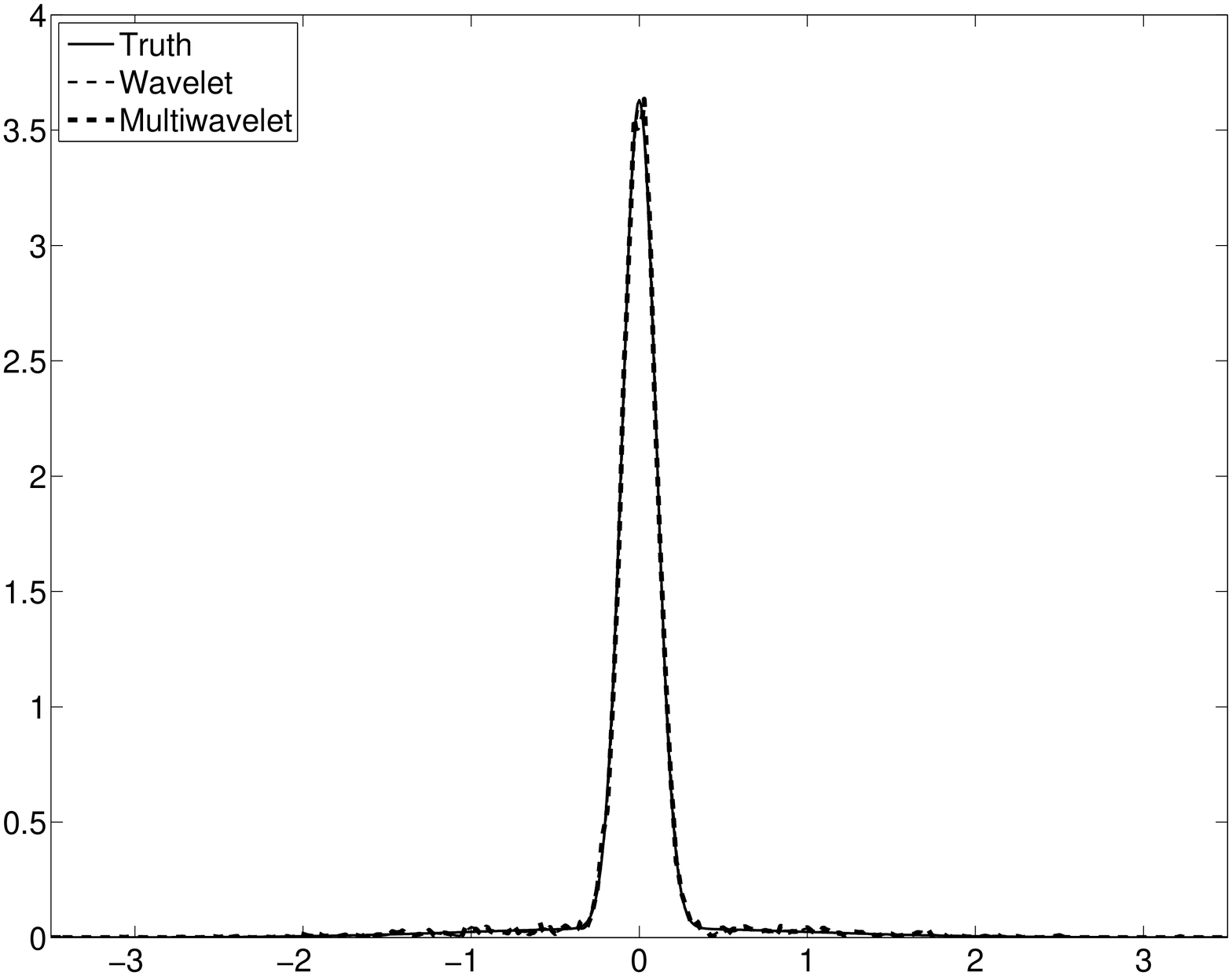}
\par\end{centering}

}\tabularnewline
\subfloat[Multisymlet 10, Symlet 4, $j=0$]{\centering{}\includegraphics[width=2in]{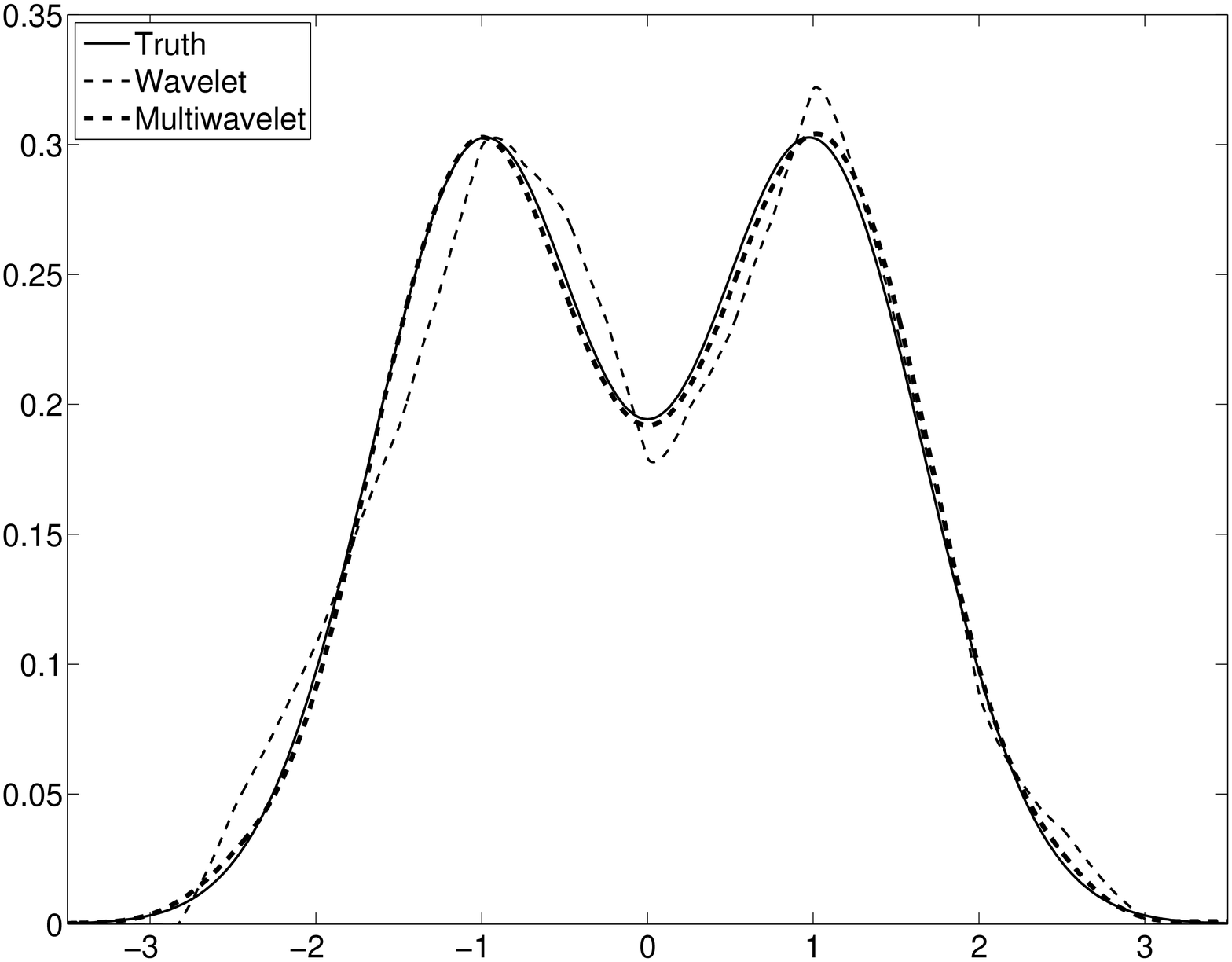}} & \subfloat[Multisymlet 7, Daubechies 4, $j=2$]{\begin{centering}
\includegraphics[width=2in]{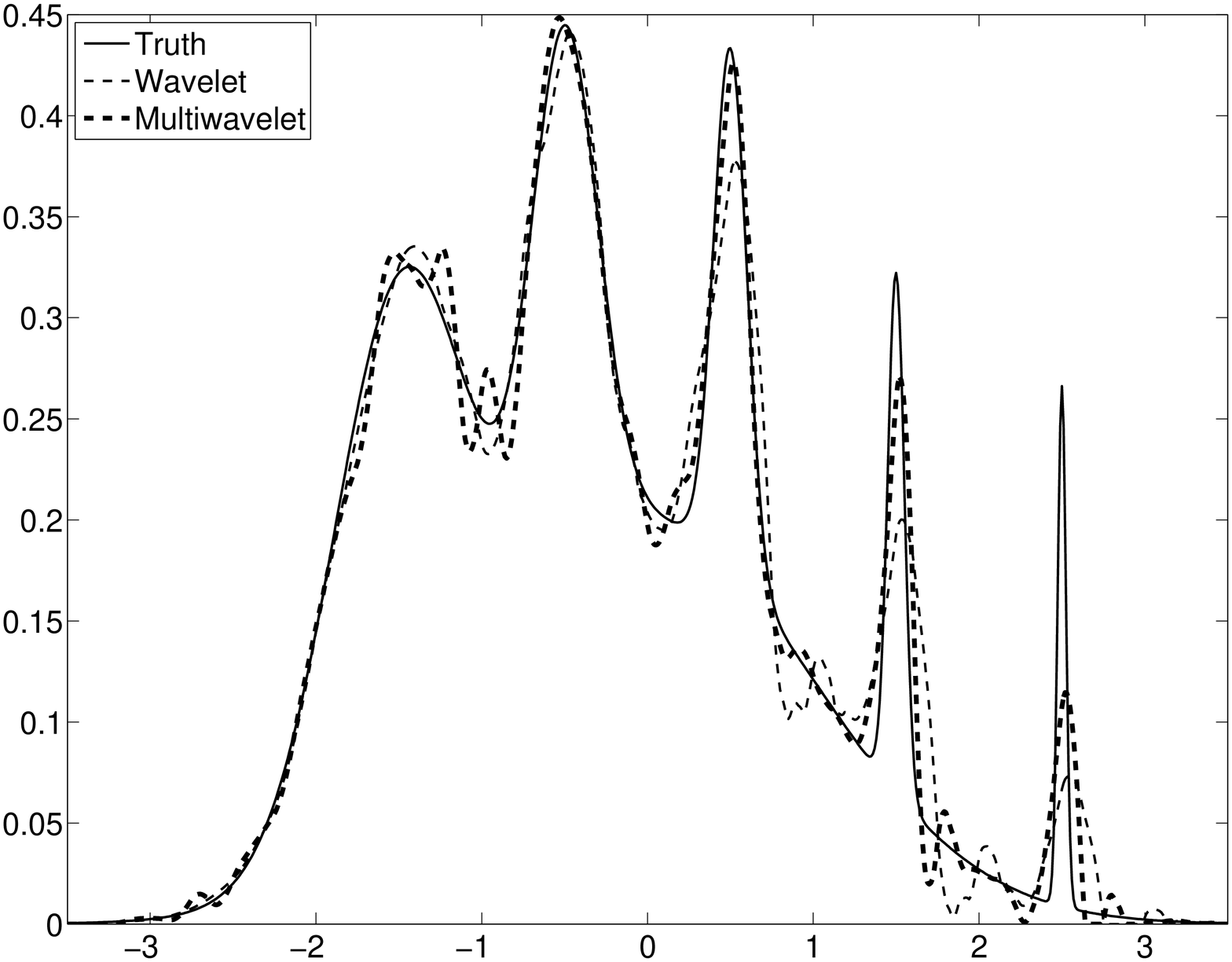}
\par\end{centering}

}\tabularnewline
\subfloat[DGHM, DMey, $j=3$]{\begin{centering}
\includegraphics[width=2in]{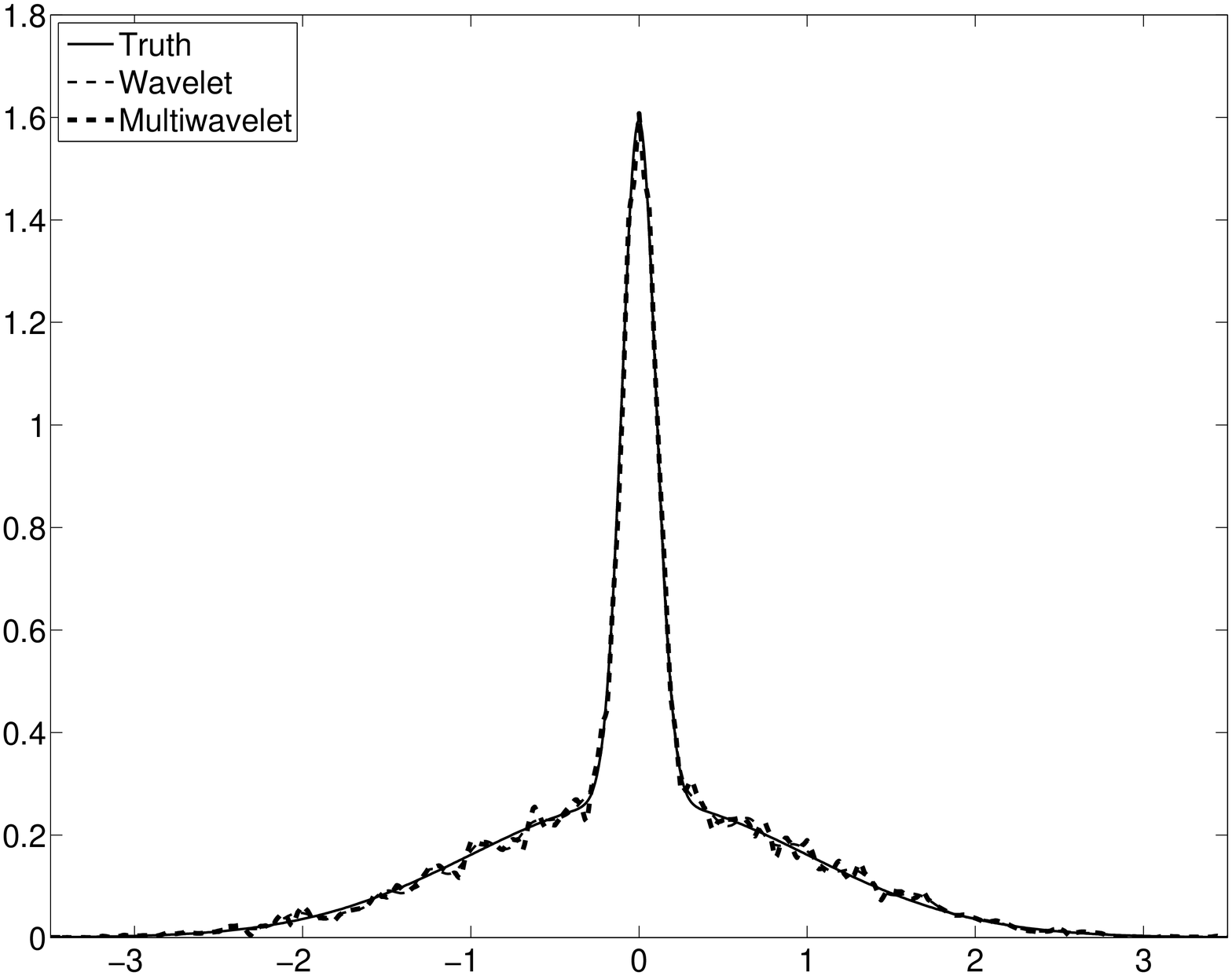}
\par\end{centering}

} & \subfloat[STT, Coiflet 5, $j=2$]{\begin{centering}
\includegraphics[width=2in]{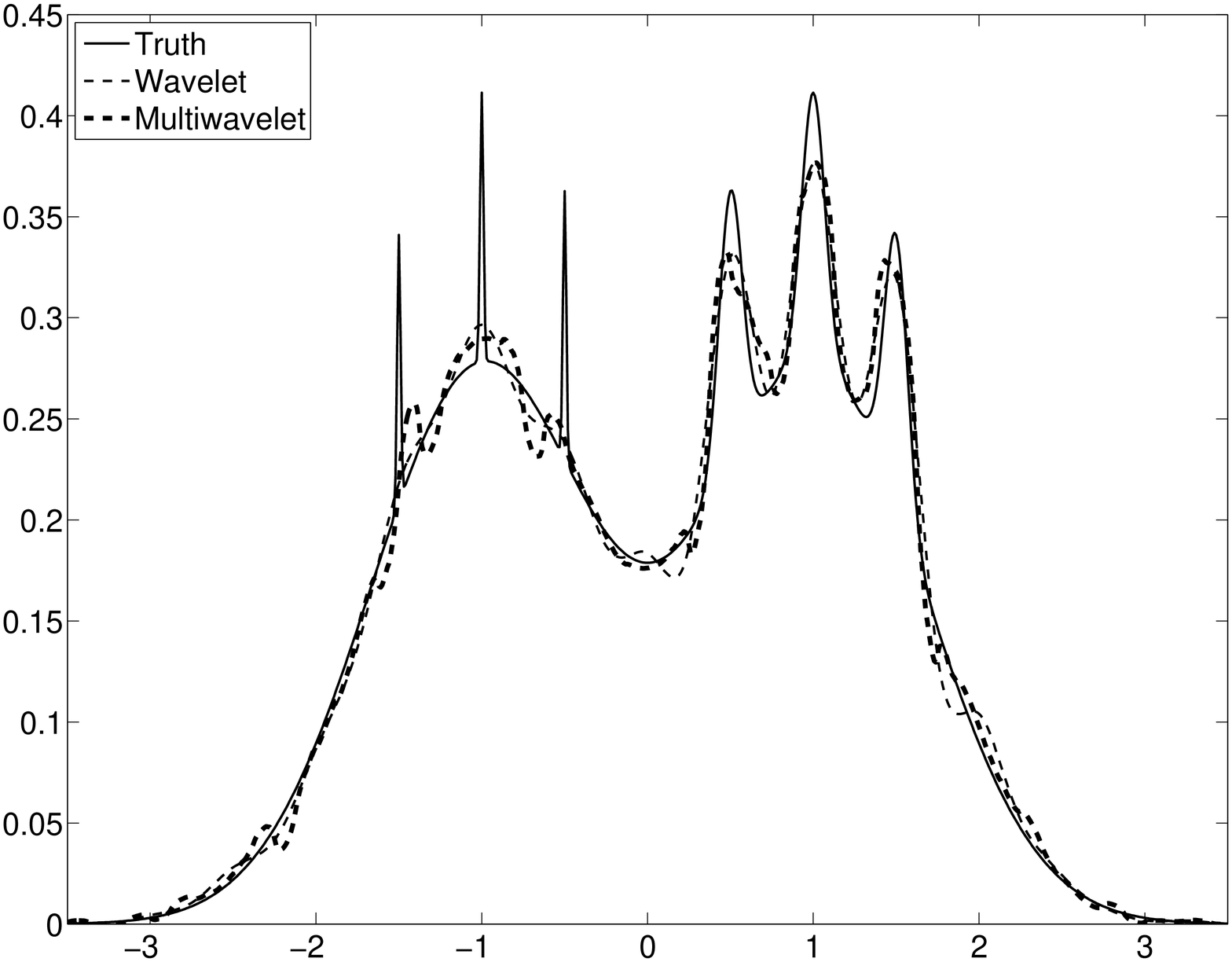}
\par\end{centering}

}\tabularnewline
\end{tabular}
\par\end{centering}

\centering{}\caption{Various densities, each with 10000 samples, approximated with a variety
of multiwavelets and wavelets to showcase the utility of MWDE alongside
the familiar WDE presented purely as a benchmark. The multiwavelet
family, wavelet family, and resolution level $j$ are provided, respectively,
in each sub-figure. \label{fig:Showcase}}
\end{figure}

There are many families of multiwavelets available, and all the orthogonal
families can be utilized for density estimation using the methods
presented in this paper. To show the utility of MWDE, we estimated
a broad range of densities from \citet{Marron92} and \citet{Wand93}
using a variety of multiwavelet bases. Along with the MWDE, we show
reconstructions using standard wavelet bases on the same plot; the
WDE are presented here purely as benchmarks, not for pitting wavelets
against multiwavelets. The results of these experiments are showcased
in Fig. \ref{fig:Showcase}. It is worth noting that the WDE and MWDE
are performed with the same resolution level in each density estimation
plot. As we have seen in the previous two sections, the most evident
difference in MWDE and WDE occurs when comparing across resolution
levels. That is, MWDE tends to perform better than WDE at coarser
resolutions.

We also conducted experiments using the cross product of all the multiwavelets
available in \citet{Keinert04} and the wavelets in \citet{PeterToolbox}
(these wavelets and multiwavelets are listed in Tab. \ref{tab:WaveletFamilies})
across resolutions $-2\leq j\leq3$ and for a variety of density functions
from \citet{Marron92} and \citet{Wand93}. The parameters of the
best (measured with $\mathrm{ISE}$) MWDE and WDE reconstructions
of each density are shown in Tab. \ref{tab:NumericalResults}. The
STT multiwavelet was the most successful multiwavelet for density
estimation, as may be expected from its smooth appearance. In addition
to the wavelet/multiwavelet families, resolution level, and $\mathrm{ISE}$
for each density, we also list the number of coefficients needed in
the expansion to achieve the given result. Note here that no MRA (thus
no thresholding) has been implemented, so the coefficient counts presented
are simply the numbers of coefficients required by the basis to span
the domain of the sample. 

\begin{table}
\noindent \begin{centering}
\begin{center}
\begin{tabular}{cc}
Wavelet Families & Multiwavelet Families \\
\hline
Daubechies 2--10 & Balanced Daubechies 2--10 \\
Symlets 4--10 & Multisymlets 4--10 \\
Coiflets 1--5 & Chui-Lian 2--3\\
Discrete Meyer & Bat \\
 & DGHM \\
 & STT  \\
\hline
\end{tabular}
\end{center}
\par\end{centering}

\caption{Wavelet and multiwavelet families used in the density estimations
in Tab. \ref{tab:NumericalResults}. \label{tab:WaveletFamilies}}
\end{table}

\begin{table}
\noindent \begin{centering}
\begin{center}
\begin{tabular}{c|ccc|ccc}
%\toprule
        &     &    MWDE    &        &     &    WDE     &        \\
Density & ISE & Res. & Coeff. & ISE & Res. & Coeff. \\
        & $(\times 10^{-3})$ & $j$ & \# & $(\times 10^{-3})$ & $j$ & \# \\
\hline
%\midrule
Normal & 0.576 & 0 & 24 & 0.194 & 0 & 38 \\
Bimodal & 0.230 & 0 & 24 & 0.124 & 1 & 44 \\
Skewed Bimodal & 0.183 & 1 & 40 & 0.0997 & 1 & 46 \\
Claw & 1.25 & 3 & 128 & 0.659 & 3 & 90 \\
Double Claw & 1.67 & 2 & 64 & 1.33 & 3 & 85 \\
%\bottomrule
\end{tabular}
\end{center}
\par\end{centering}

\caption{Results of estimating a variety of densities with both MWDE and WDE
with several families of wavelets across resolutions $-2\leq j\leq3$.
Information about the density estimation with the lowest $\mathrm{ISE}$
is shown for both wavelets and multiwavelets. The wavelet and multiwavelet
families used are listed in Tab. \ref{tab:WaveletFamilies}.\label{tab:NumericalResults}}
\end{table}

\section{Discussion\label{sec:Discussion}}

We are motivated to investigate MWDE, for we can use symmetric and
symmetric/antisymmetric bases, such as the STT and DGHM multiwavelets,
instead of being limited to the necessarily asymmetric orthogonal
wavelets. So, we began our investigation by comparing the results
of using the DGHM multiwavelet for MWDE and the Daubechies wavelet
of order 2 for WDE for estimating a symmetrical Gaussian distribution.
As evident from Fig. \ref{fig:DGHMandGauss}, the DGHM multiwavelet
outperforms the Daubechies wavelet at coarser resolution levels. Following
this, we used the more suitable STT multiwavelet and Daubechies wavelet
of order 5 to estimate the same density with much better results.
The STT multiwavelet is much smoother than the DGHM multiwavelet,
as is evident from comparing Figs. \ref{fig:DGHM} and \ref{fig:STT}.
We see in Fig. \ref{fig:STTandGauss} that both the MWDE and WDE perform
well, particularly at resolutions $j=0$ and $j=1$, respectively.
The primary point of this empirical result is that the MWDE performed
its best at a coarser resolution level than the wavelet. We will find
this is a recurring theme throughout the tests we performed. This
is an interesting relationship, one which we will continue to explore
and develop more formally in later works. 

We also investigated a specific and interesting class of multiwavelets:
the so-called ``balanced multiwavelets.'' We tested MWDE using balanced
multiwavelet bases and display the results in Fig. \ref{fig:BalProgression}.
The objective of this exercise was to examine the relationship across
resolution levels between MWDE and WDE using related bases. As shown
by \citet{Lebrun1998}, we can balance a standard Daubechies wavelet
to a balanced Daubechies multiwavelet. We did this using the Daubechies
wavelet of order 5 balanced to a multiwavelet of multiplicity 2. We
then compared across resolutions the density estimations resulting
from using the Daubechies wavelet of order 5 for WDE and the balanced
Daubechies multiwavelet of order 5 for MWDE. We found that MWDE and
WDE result in very similar density estimates, but the wavelet ``lags''
the multiwavelet. That is, the MWDE at resolution level $j$ is very
close to the WDE at resolution $j+1$. This is expected in the case
of Daubechies wavelets because the balanced Daubechies multiwavelets
are just compressed and translated on the half-integers versions of
the the base Daubechies wavelet. We demonstrated this in Fig. \ref{fig:BalDaub2}.

There are, of course, many wavelet and multiwavelet families available
with the necessary properties---namely orthogonality---for density
estimation using the methods presented in this paper. We performed
MWDE and WDE using a variety of multiwavelet and wavelet families
on a broad sample of distributions and across several resolution levels.
Our results are summarized in Tabs. \ref{tab:WaveletFamilies} and
\ref{tab:NumericalResults}. We consistently see the STT multiwavelet
and Coiflet 5 wavelet perform the best under the $\mathrm{ISE}$.
What is perhaps more interesting are the numbers of coefficients required
in the various expansions to produce the density estimations. These
are explicitly given in Tab. \ref{tab:NumericalResults}. As we have
not used any non-linear density estimation, the numbers of coefficients
presented are simply the numbers of coefficients required for the
scaling and multiscaling functions to span the sample (which, in most
cases, is contained on the unitless interval $\left[-4,4\right]$).
From the table, we see that WDE, under the $\mathrm{ISE}$, provides
the best density estimation for the given distributions. However,
in most cases, MWDE provides its best results at a coarser or equal
resolution and, in all but one case, using fewer coefficients than
the best WDE results. This could prove fruitful in terms of sparse
representation if an MRA is constructed and non-linear density estimation
is performed using multiwavelet bases. Finally, in Fig. \ref{fig:Showcase},
in order to show the utility of MWDE, we show a multitude of density
estimations on various densities and using various multiwavelet families.

In the cases we have investigated, we see that multiwavelets provide
their best density estimation at resolution levels coarser than or
equal to the best wavelet estimation for any particular density. This
is expected, as a multiwavelet density estimator is constructed such
that there are multiple basis functions at every translate along the
domain requiring two (or, generally, $r$) coefficients for every
translate instead of just one as is needed by wavelets. In summary,
our investigation shows a general trend that the MWDE converges to
its best estimate at resolution levels coarser than or equal to comparable
WDE and using a comparable or fewer number of coefficients.

\section{Conclusion\label{sec:Conclusion}}

In this paper, we have presented for the first time the use of multiwavelet
bases for density estimation. We illustrated how to implement multiwavelet
density estimation (MWDE) by projecting onto orthogonal multiwavelet
bases. The utility of the approach was demonstrated by estimating
several densities and comparisons were conducted with the well-established
wavelet density estimation (WDE) as a benchmark. Our results show,
in the large, that MWDE converges to its best estimate at resolution
levels coarser than or equal to comparable WDE. Furthermore, the number
of coefficients required by the best MWDE was, in all but one case,
less than the number of coefficients required by the best WDE. We
also examined MWDE using balanced multiwavelets and made some interesting
empirical observations. We showed that WDE ``lagged'' MWDE by one
resolution level when a balanced Daubechies multiwavelet of multiplicity
2 is used for the MWDE and the corresponding Daubechies wavelet is
used for the WDE. Furthermore, we showed the STT multiwavelet was
the best (measured by $\mathrm{ISE}$) basis of the families we tested
for estimating a broad range of densities. In future research, we
plan to investigate non-linear MWDE through incorporation of vector
thresholding techniques and direct non-negative density estimation
by estimating $\sqrt{p}$ expanded in a multiwavelet basis.

\section{Acknowledgments}

The authors thank Dr. Fritz Keinert of Iowa State University for illuminating
communications.

This research is based on work supported by the National Science Foundation's
Research Experience for Undergraduates program under grant IIS-REU-0647018.
Any opinions, findings, and conclusions or recommendations expressed
in this material are those of the author(s) and do not necessarily
reflect the views of the National Science Foundation.

\section{References}

\bibliographystyle{elsarticle-harv}
\bibliography{references}

\end{document}